\documentclass[a4paper,11pt]{article}   

\usepackage{amsmath,amsfonts,amssymb,amsthm,amsxtra}
\usepackage{array}

\newcount\pictures
\ifnum\pictures=1
\usepackage{epsfig}
\newcommand{\wherepoly}{ps_files}
\fi
\newcommand{\inte}{{\rm int\,}}
\newcommand{\bd}{{\rm bd\,}}

\newcommand{\conv}{{\rm conv\,}}
\newcommand{\aff}{{\rm aff\,}}
\newcommand{\lin}{{\rm lin\,}}

\newcommand{\res}{{\rm res\,}}

\newcommand{\Kd}{{\cal K}^d}

\newcommand{\R}{{\mathbb R}}
\newcommand{\N}{{\mathbb N}}

\newcommand{\Z}{{\mathbb Z}}
\newcommand{\C}{{\mathbb C}}

\newtheorem{theorem}{Theorem}[section]
\newtheorem{lemma}{Lemma}[section]

\newtheorem{remark}{Remark}[section]

\newtheorem{corollary}{Corollary}[section]

\newtheorem{algorithm}{Algorithm}[section]
\numberwithin{equation}{section}

\begin{document}

\title{Densest Lattice Packings of 3--Polytopes}
\author{Ulrich Betke \and Martin Henk\thanks{The second named author acknowledges the hospitality of the International Erwin Schr\"odinger  Institute for Mathematical Physics in Vienna, where a main part of his contribution to  this work has been completed.}}

\date{}
\maketitle

\begin{abstract}
Based on Minkowski's work on critical lattices of 3-dimensional convex bodies we
 present an efficient algorithm for computing the density of a densest lattice packing of 
an arbitrary 3-polytope. As an application we calculate densest lattice packings of all regular and  Archimedean polytopes.
\end{abstract}

\section{Introduction}\label{sec:Introduction}
Throughout this paper, $\R^d$ denotes the $d$--dimensional
Euclidean space with origin $0$, Euclidean norm $\|\cdot\|$, inner
product $\langle \cdot,\cdot \rangle$ and unit sphere
$S^{d-1}$. $\Kd$ denotes the set of all convex bodies $K\subset \R^d$ with nonempty interior $\inte(K)$ and $\Kd_0$ denotes the subset
of $\Kd$ consisting of all bodies which are centrally
symmetric with respect to the origin. For a set $M \subset \R^d$ we denote
by ${\rm vol}(M)$ its volume with respect to its affine hull $\aff(M)$. Furthermore, $\conv(M)$, $\lin(M)$ denotes the convex hull, linear hull  of $M$, respectively. The boundary of $K\in\Kd$ is denoted by $\bd(K)$.

By a lattice $\Lambda\subset \R^d$ with basis $B=\{b^1,\dots,b^d\}$, where $b^1,\dots,b^d\in \R^d$ are linearly independent,  we understand the  set
$$
\Lambda = \left\{ z_1b^1+\dots z_d b^d : 
z_1,\dots,z_d \in \Z\right\}=B\Z^d.
$$
The determinant $\det \Lambda$ of $\Lambda$ is the volume of the
parallelepiped spanned by $b^1,\dots,b^d$, i.e., $\det\Lambda=|\det B|$.  
A lattice $\Lambda$ is called a
packing lattice for $K \in \Kd$ if $x+K$ and $y+K$ do not
overlap for $x,y \in \Lambda$, $x\ne y$. There always exists a
packing lattice $\Lambda^*(K)$ with minimal determinant. Such a lattice
 is called a densest packing lattice and the
quantity
$$
\delta^*(K) = {\rm vol}(K) / \det \Lambda^*(K)
$$
is called the lattice packing density of $K$ or the density of a densest lattice packing of $K$. There is a
large amount of literature on packings. Some books, which
give a good description of the background of our work, are \cite{ErdoesGruberHammer}, \cite{GruberLekkerker}, \cite{Rogers:packing}, \cite{Zong:strange} as well as the Diploma thesis \cite{Haus:Minkowski}. 
For a more general survey on the theory of packings we refer to \cite{fejestothkuperberg:packing} and the references within.

The lattice packing problem for a general body in $\Kd$ is
very hard. In fact, for $d\ge 4$ the only exact results are
on space fillers (cf.~\cite{McMullen:spacefiller}, \cite{Venkov:spacefiller}) for which 
$\delta^*(K)=1$ and on the unit ball $B^d$ where $\delta^*(B^d)$ is
known for $d \le 8$ (see e.g.~\cite{ConwaySloane:spheres}, \cite{Zong:spheres}).
In contrast for a fixed body $K\in {\cal K}^2$ there are several techniques to solve the problem (see e.g. \cite{GruberLekkerker}, pp.~241) 
and there exists also an algorithm, due to {\rm Mount {\rm and} Silverman} \cite{MountSilverman}, that determines $\delta^*(P)$ for a centrally symmetric $n$-gon in time $O(n)$. 
   
However, already in 3--space the situation is rather more complicated. Apart from the 3-dimensional space-fillers (for a classification see \cite{GruberLekkerker}, pp. 164) and from cylinders based on a convex disk, in which case it can be shown that  the problem is equivalent to the determination of the lattice packing density of the convex disk (cf.~\cite{Rogers:packing}, p.~13), densest lattice packings  are  only known for the  bodies in $\R^3$, which are listed in the following table:

\vspace{0.3cm}
\begin{center}
\begin{tabular}{ m{0.54\hsize} | >{\quad} m{0.37\hsize} }
  Body & Author \\ \hline 
  regular octahedron & {\rm Minkowski} \cite{m}, 1904 \\[0.2ex]
  truncated cubes; for $0< \lambda \leq 3$ : \hfill \allowbreak
  $\left\{ x \in \R^3 : |x_i|\le 1,\; |x_1+x_2+x_3|\le \lambda\right\}$  &  {\rm Whitworth} \cite{w1}, 1948 \\[0.2ex]\vspace{0.1cm}
$\left\{ x\in \R^3 : \left(x_1^2+x_2^2\right)^{1/2} +|x_3| \le 1\right\}$ & {\rm Whitworth} \cite{w2}, 1951 \\[0.2ex]\vspace{0.1cm}
frustrum of a sphere; for $0<\lambda\leq 1$: \hfill \allowbreak 
 $\left\{ x \in \R^3 : x_1^2+x_2^2+x_3^2 \le
1,\; |x_3| \le \lambda \right\}$ & {\rm Chalk} \cite{c}, 1950 \\[0.5ex]\vspace{.1cm}
tetrahedron and cubeoctahedron & {\rm Hoylman}  \cite{h}, 1970 
\end{tabular}
\end{center}

It is worth to mention that the family of frustrums of a sphere includes the 3-ball as a limiting case,
for which the packing density was determined already by {\rm Gauss} \cite{g}. 

All the computations of $\delta^*(K)$ for the bodies in the table above may be regarded as an application of a general method developed by {\rm Minkowski} \cite{m} which characterizes densest packing lattices of a 3-dimensional convex body by certain properties (see also \cite{GruberLekkerker}, pp.~340).  However this  method was considered
as rather impractical and in 1964 {\rm Rogers} \cite{Rogers:packing} wrote ``{\em Despite considerable theoretical advances in 
  the Geometry of Numbers since Minkowski's time, the
  problem of determining the value of $\delta^*(K)$ for a
  given convex $3$-dimensional body $K$ remains a formidable task}''. Indeed the only change in the list of known densest lattice packings since the publication of {\rm Rogers'} book is the addition of the tetrahedron and  cubeoctahedron. And even in 1990 {\rm Gruber} mentioned the determination of $\delta^*(K)$ of a $3$-dimensional convex body as one of the important open problems in Geometry of Numbers \cite{Gruber:dmvarticle}.

Here we use {\rm Minkowski}'s work as a starting point for the construction of a
practicable algorithm to compute the packing density of an
arbitrary polytope in $\R^3$. We proceed as follows. In
section 2 we adapt {\rm Minkowski}'s work to our purposes.
Moreover we discuss in this section the principal
applicability of his work in higher dimensions.

Having adapted {\rm Minkowski}'s method we still have to solve two problems to obtain a practicable algorithm. First we have to determine
the local minima of a polynomial of degree three in
three variables. Somewhat surprisingly there appears to be
no general purpose algorithm of numerical analysis which
suits our needs. Thus we develop an ad hoc method for our
problem. This will be done in section 4. Before doing this
we deal in section 3 with a more geometric problem. It
turns out that by the methods of section 2 we have
essentially to look at every choice of 7 facets of the
polytope. This leaves us with more or less $\binom{n}{7}$ cases for a polytope
with $n$ facets and in every step we have to solve a
nonlinear minimization problem. This  limits the
applicability of the algorithm to  polytopes with   few facets. Thus
we develop in section 3 some methods to reduce the number
of cases. While we do not give an exact worst case analysis
it should be possible to reduce the number of cases to
possibly as few as $O(n^{3/2})$ in typical cases  by showing
that optimal packings must lead to feasible points of some
 related simple linear optimization problems. This will be
done in section 3. It should however be mentioned that in
our program we do not fully exploit this reduction as an
implementation would become  rather complex. As an application of our work we present in
section 5 optimal lattice packings for all
regular and Archimedean polytopes.

\section{Necessary conditions for optimal lattices}\label{section2}

In this section we state without proof {\rm Minkowski}'s results.  In
fact, we have summarized the relevant facts in Theorem \ref{ness_con}.
We have included some more results, as we discuss at the end of the
section, whether the algorithm could be extended to higher dimensions.
Further we have changed and modernized the notation. We remark that a
good part of {\rm Minkowski}´s theory is exposed in \cite{GruberLekkerker} and
\cite{Zong:strange}, though in both books Theorem \ref{ness_con} is stated in a slightly  weaker form. 

For $K_i\in{\cal K}^d$, $\lambda_i\in\R$, $i=1,2$, we
denote by $\lambda_1 K_1+\lambda_2 K_2$ the set 
$ \{\lambda_1 x^1 +\lambda_2 x^2 : \allowbreak 
x^1 \in K_1, \; x^2 \in K_2\}$. {\rm Minkowski} observed, that 
\begin{equation}
\label{sym_body}
\Lambda \text{ \it is a packing lattice of } K \Leftrightarrow 
 \Lambda \text{ \it is a packing lattice of } \frac{1}{2}(K  - K).
\end{equation}
Since the difference body $\frac{1}{2}(K  - K)$ belongs to the class ${\cal K}_0^d$, in the following we assume that all bodies are centrally symmetric.

A lattice $\Lambda$ is called admissible for 
$K\in \Kd_0$, if $\inte K \cap \Lambda=\{0\}$. It is well
known and easy to see that $\Lambda$ is admissible if and only if $2
\Lambda$ is a packing lattice. The value
$$
\Delta(K) = \min\{\det \Lambda : \Lambda \text{ admissible for
  }K\}
$$
is called the critical determinant of $K$ and an admissible lattice
$\Lambda$ satisfying $\Delta(K)= \det \Lambda$ is called a critical
lattice. Thus $\delta^*(K)=\text{\rm vol}(K)/(2^d\Delta(K))$ and the problems of constructing a 
densest packing lattice and a critical lattice are equivalent.

While we naturally do not know a basis of a critical lattice
beforehand, we shall see that we have a great amount of information on
the behaviour of certain lattice points with fixed coordinates with
respect to such a basis. For a basis $B=\{b^1,\dots,b^d\}$ of
$\R^d$ and a point $x\in \R^d$ we denote by $x_B=(x_1,\dots,x_d)_B^\intercal$ the vector given by $x_B=\sum_{i=1}^n x_i b^i$.

The construction of critical lattices is based on the connection of
lattice crosspolytopes and primitive vectors. As usual, a set of lattice vectors $b^1,\dots, b^k$ of a lattice $\Lambda\subset \R^d$ is called {\em primitive} iff  this set can be extended to a basis of the lattice $\Lambda$.
Let $\Lambda\subset \R^d$ be a lattice and let $b^1,\dots,b^k\in\Lambda$ be linearly independent. The crosspolytope $C=\conv\{\pm b^1,\dots,\pm b^k\}$ is called a {\em lattice crosspolytope} iff $\Lambda\cap \inte(C)=\{0\}$. If we even have  $\Lambda\cap \bd(C)=\{\pm b^1,\dots,\pm b^k\}$ then it is called a {\em free lattice crosspolytope}.
  Clearly, the convex hull of  every primitive set $b^1,\dots,b^k$
forms a free lattice crosspolytope  and every free lattice crosspolytope is a
lattice crosspolytope, while the converse is not true. A complete
characterization of (free) lattice crosspolytopes of dimension up to three
is given in

\begin{lemma}[Minkowski, 1904]
  Let $\Lambda$ be a lattice in $\R^3$.
\begin{itemize}
\item  For $k=1$ the vertices of any lattice crosspolytope are primitive.
\item   For $k=2$ the vertices of any free lattice crosspolytope  are primitive,
  while for the non--free lattice crosspolytopes $C$ there exists a basis $B$
  of the lattice such that $C=\conv\{(1,0,0)_B, (1,2,0)_B\}$.
\item  For $k=3$ there are two types of free lattice crosspolytopes. One has
  primitive vertices {\em (crosspolytope of the first type)}, the other has the
  vertices $(1,0,0)_B,(0,1,0)_B,(1,1,2)_B$ for a basis
  $B$ {\em (crosspolytope of the second type)}. For every non--free lattice
   crosspolytope $C$ there exists a basis $B$ such that
  $C=\conv\{(1,0,0)_B,(0,1,0)_B,p\}$, where $p$ is element of  the set 
  $\big\{(0,1,2)_B,(1,1,2)_B,(1,2,2)_B,(1,1,3)_B, (1,2,3)_B$, \linebreak[5]  $(2,2,3)_B$,
     $(1,2,4)_B$,  $(2,3,4)_B\big\}$.
\end{itemize}
\end{lemma}
   
\noindent Using free lattice crosspolytopes we can identify critical lattices for
every $K \in {\cal K}_0^3$. To this end we use the following abbreviation:  
 For a basis  $B=\{b^1,b^2,b^3\}$ of $\R^3$ let 
 
\begin{equation} 
\label{def_testsets}
\begin{split}
  {\cal U}_B^1 & = \big\{(1,0,0)_B,\;\allowbreak (0,1,0)_B,\;(0,0,1)_B,\;
        (0,1,-1)_B, (-1,0,1)_B, (1,-1,0)_B\big\}, \\
 {\cal U}_B^2 & = \big\{(1,0,0)_B,\;(0,1,0)_B,\;(0,0,1)_B,\;(0,1,1)_B,\; 
  (1,0,1)_B,\; (1,1,0)_B\big\}, \\
  {\cal U}_B^3 & =\big \{(1,0,0)_B,\;(0,1,0)_B,\;(0,0,1)_B,\;(0,1,1)_B,\;
  (1,0,1)_B,\;(1,1,0)_B,\\  &\quad\quad (1,1,1)_B\big\}.
\end{split}
\end{equation}

\begin{lemma}[Minkowski, 1904]
\label{pack_char}
  Let $K \in {\cal K}^3_0$. Then there exists a critical lattice
  $\Lambda$ with basis $B$ such that one of the following cases holds:
\begin{enumerate}
\item[{\rm 1.}] ${\cal U}_B^1 \subset \bd(K)$ and $K$ contains no lattice crosspolytope of the second type.
\item[{\rm 2.}]
  ${\cal U}_B^2\subset \bd(K)$ and $(1,1,1)_B \not \in K$.
\item[{\rm 3.}]
 $ {\cal U}_B^3\subset \bd(K)$.
\end{enumerate}
\end{lemma}

\noindent While the distinction between cases 2. and 3. may look
artificial, we shall see that it leads to rather different cases in
the actual computation of a critical lattice.
For any lattice with basis $B$ which satisfies a condition of the previous lemma we can easily check its admissibility:

\begin{lemma}[Minkowski, 1904] 
\label{lattice_check}
Let $K \in {\cal K}_0^3$ and $\Lambda$ be a lattice with basis $B$. Then $\Lambda$ is admissible, if one of the following conditions is satisfied
\begin{enumerate} 
\item[{\rm 1.}] 
  ${\cal U}_B^1 \subset \bd(K)$ {\em and} 
  $(-1,1,1)_B, (1,-1,1)_B, (1,1,-1)_B \not\in  \inte K$. 
\item[{\rm 2.}] ${\cal U}_B^2 \subset \bd(K)$ {\em and}  $(1,1,1)_B \not\in \inte K$.
\end{enumerate}
\end{lemma}
\begin{proof} The second statement follows immediately from Minkowski's characterizations of lattice crosspolytopes  (cf.~\cite{Zong:strange}, Lemma 4.9, $2^*$), whereas item  {\rm 1.~}is contained in a more implicitly way in his paper.   
However, from the work of Minkowski we know  (cf.~\cite{Zong:strange}, Lemma 4.9, $1^*$) that in the first case, i.e., ${\cal U}_B^1 \subset \bd(K)$, the lattice $\Lambda$ is admissible if 
\begin{equation}
\label{case1}
  \text{\rm a1) }\, (0,\pm1,\pm1)_B \notin K,\quad 
  \text{\rm a2) }\, (1,\pm1,\pm1)_B \notin K,\quad 
  \text{\rm a3) }\, (2,\pm1,\pm1)_B \notin K, 
\end{equation}
as well as all points arising from permutations of the coordinates of the above points are not contained in $K$.
Therefore, in order to prove {\rm 1.} we have to show that the points $(0,1,1)_B$, $(1,1,1)_B$, $(2,\pm 1, \pm 1)_B$ and the corresponding permutations are not contained in $K$. To this end let  $f: \R^3\to \R_{\geq 0}$ be the distance function of $K$, i.e, $f(x)=\min\{ \lambda : \lambda \in \R_{\geq 0} \text{ and } x \in \lambda K\}$. Then $K=\{ x \in \R^3 : f(x)\leq 1\}$ and since $K\in {\cal K}_0^3$ the function $f$ describes a norm on $\R^3$. Since ${\cal U}_B^1 \subset \bd(K)$ we find 
\begin{align*}
 f((0,1,1)_B)\geq f((0,0,2)_B)-f((0,-1,1)_B) & =1,\\
 f((1,1,1)_B)\geq f((3,0,0)_B)-f((1,-1,0)_B)-f((1,0,-1)_B) & =1, \\
 f((2,1,-1)_B) \geq f((2,0,0)_B)-f((0,-1,1)_B) & =1, \\
 f((2,-1,1)_B) \geq f((2,0,0)_B)-f((0,1,-1)_B) & =1, \\
 f((2,-1,-1)_B) \geq f((2,-2,0)_B)-f((0,-1,1)_B)) & =1.  
\end{align*} 
Obviously, the same inequalities hold for the points given by all permutations of the coordinates of the points of the left hand side and thus 
\begin{equation*}
 f((2,1,1)_B) \geq f((2,2,0)_B)-f((0,1,-1)_B)  \geq 1.
\end{equation*} 
\end{proof}

\noindent We  call the sets ${\cal U}_B^j$, $j=1,2,3$, (cf.~\eqref{def_testsets}) {\em test sets} of the first, second or  third kind, respectively.   

Now suppose that $B=\{b^1,b^2,b^3\}$ is a basis of a critical lattice $\Lambda$ of $K$ and  let ${\cal U}_B=\{u^1_B,\dots,u^k_B\}$, $k=6$ or $k=7$, be one the three test sets such that ${\cal U}_B\subset \bd(K)$ (cf.~Lemma \ref{pack_char}). Let $H_i$ be {\em any} supporting hyperplane of $K$ containing $u^i_B$, $i=1,\dots, k$, and let     
\begin{equation}
\label{set_S}
  {\cal S}_{H_1,\dots, H_k}=\left\{ W\in \R^{3\times 3} : u^i_W \in H_i, \; 
                                   1\leq i\leq k\right\}. 
\end{equation}
On this space we consider the function  
\begin{equation}
\label{function}
\begin{split}
   f_{H_1,\dots,H_k}: {\cal S}_{H_1,\dots, H_k} &\to \R \text{ given by } \\ 
   & f_{H_1,\dots,H_k}(W)=|\det(W)|.
\end{split}
\end{equation}
If we  assume for a moment that ${\cal U}_B=\Lambda\cap\bd(K)$ then 
\begin{equation}
\label{loc_min}
  B \text{ is a local minimum  of    }  f_{H_1,\dots,H_k}(W), \quad W\in {\cal S}_{H_1,\dots, H_k}.
\end{equation} 
Otherwise there exists a $W\in {\cal S}_{H_1,\dots, H_k}$ in a sufficiently small neighbourhood of $B$ such that $|\det(W)|<\det\Lambda$ and $K \cap (W\Z^d)\backslash\{0\}\subset\{u^1_W,\dots, \allowbreak  u^k_W\}$. However, this set of vectors is contained in the supporting hyperplanes $H_1,\dots,H_k$ for any  element $W\in {\cal S}(H_1,\dots, H_k)$. Thus the lattice $W\Z^d$ is admissible and hence $\Lambda$ can not be critical.

In general the situation is much more complicated as ${\cal U}_B\cap
\bd(K)$ may just be a proper subset of $\Lambda \cap \bd(K)$. 
But by a
close examination of all possible cases {\rm Minkowski} found

\begin{theorem}[Minkowski, 1904]
 \label{ness_con}
  Let $K \in {\cal K}_0^3$. Then there exists a critical lattice
  $\Lambda$ of $K$ with basis $B$ such that $\bd(K) \cap \Lambda$ contains a test set ${\cal U}^j_B=\{u^1_B,\dots,u^k_B\}$ for a   $j\in\{1,2,3\}$, such that for any choice of supporting hyperplanes $H_i$ of $K$ containing $u_B^i$, $1\leq i \leq k$,  one of the following 4 cases holds:
\begin{enumerate}
\item [{\rm I.}]  $j=1$  and \eqref{loc_min} holds.
\item [{\rm II.}]  $j=2$  and \eqref{loc_min} holds.
\item [{\rm III.}]  $j=3$  and \eqref{loc_min} holds.
\item [{\rm IV.}]  $j=3$  and there are scalars 
                         $\lambda_1,\lambda_2,\lambda_6>0$, $\lambda_3\in\R$, such that   for the outer normal vectors $v^1,v^2,v^3,v^6$ of the hyperplanes $H_1,H_2, H_3,H_6$, respectively, holds: {\rm a)} $\lambda_1 v^{1}+\lambda_2 v^{2}+\lambda_3 v^{3}=\lambda_6v^{6}$ and  {\rm b)} the hyperplane $\tilde{H}_6$ with outer normal vector $\lambda_1 v^{1}+\lambda_2 v^{2}$ containing $u_B^6$ is a supporting hyperplane of $K$  and \eqref{loc_min} holds with $H_6$ replaced by ${\tilde H}_6$.
\end{enumerate}
\end{theorem}

At this point we should remark that in  case I.~{\rm Minkowski}
gave a somewhat stronger condition, but our form appears to be
more suitable for automatic computation.

Of course, given an arbitrary convex body $K$ we  do not know how to exploit Theorem \ref{ness_con}, but if we consider only polytopes  then  
for the  supporting hyperplanes $H_i$ in Theorem \ref{ness_con} we may always choose the supporting hyperplanes of the facets of the polytope. As a polytope has only finitely many 
facets we obtain the following frame of an algorithm for the computation 
of a critical lattice of a polytope:

\begin{algorithm} Let $P\in{\mathcal K}^3_0$.\hfill
\label{algorithm_org}
\vspace{-0.2cm}
\begin{itemize}
\item For each of the  cases {\rm I., II., III., IV.}~of Theorem \ref{ness_con} do \vspace{-0.2cm}
\begin{itemize}
\item For every choice of $k$ facets with supporting planes $H_1,\dots,H_k$ of $P$  {\rm (}$k=6$ in the first two cases
 and $k=7$ in the latter ones{\rm )}  do 
\begin{enumerate}
    \item[{\rm S1.}] Determine ${\cal S}_{H_1,\dots, H_k}$.
    \item[{\rm S2.}] Find the local minima ${\cal M}_{H_1,\dots,H_k}$ of $f_{H_1,\dots,H_k}$ {\rm (}cf.~\eqref{loc_min}{\rm )}. 
    \item[{\rm S3.}] For each  $M\in {\cal M}_{H_1,\dots,H_k}$ check whether $M\cdot \Z^d$ is an admissible lattice, i.e., $M$ satisfies the criterion {\rm 1}.~of Lemma {\rm \ref{lattice_check}} in the first case and criterion {\rm 2.}~in the remaining cases. 
\end{enumerate}
\end{itemize}
\item Among all calculated admissible lattices find one with  minimal determinant. The corresponding lattice is a critical lattice of $P$.
\end{itemize}
\end{algorithm}

It turns out that at this point we are left with two problems. First
there appears to be no general purpose algorithm to find all local
minima of a function  like $f_{H_1,\dots,H_k}$ as we have to do in Step S2. Moreover, a priori we can not assume that the local minima of this function are isolated points. In general, they may form a manifold and  we have the problem to parameterize such a  manifold in order to carry out step S3.  In section \ref{min_pol} we shall show how one can overcome these problems. 
Another  problem is just the number of steps of the algorithm. In a straightforward implementation, we have to consider every choice of $6$ and $7$ facets of the polytope and hence we have about  $\binom{n}{7}$ steps.  Of course, this limits the algorithm to polytopes having only  few facets. Hence in order to get an efficient algorithm we have  to reduce the
number of steps. This can be done very effectively by some
considerations given in section \ref{num_red}.

We close the section with some remarks on the extension of the
algorithm to higher dimensions. While {\rm Minkowski} settles his work
in 3-space, the ideas principally work in higher dimensions as
well. In fact there has been an enumeration of lattice crosspolytopes  in
dimension  4 (\cite{br}, \cite{laub:gitterokateder}, \cite{wo}, \cite{mo},
\cite{ba1}, \cite{ba2}) and dimension 5 (\cite{Wessels:crosspolytopes}). 
 Beside this classification we have to determine the number of different 
test sets. For the cardinality of a test set we have the general natural lower
bound of $d(d+1)/2$ given in \cite{Swinnerton-Deyer}. There is no obvious upper bound
for their cardinality, but the results in dimension two and three
suggest the upper bound of one half of the maximal number of lattice points contained in the boundary
of a lattice point free strictly convex set. Due to a result of Minkowski this number is bounded by $2^{d}-1$.
Finally,  we have to take into account the additional lattice points (not contained in the test sets) lying in the
boundary, which are responsible for  the split of test sets of the
third kind into two separate cases in dimension 3. Thus even in dimension 4 there
should be a large number of cases which have to be considered
separately and there has been no attempt to give an enumeration of
these cases.

Further in each case we have to consider all possible choices of
facets and even though we could apply the results of the next section
it appears to be computationally difficult to choose at least $10$ out
of $n$ facets. Moreover, even in dimension 4  we have to solve a minimization problem for
polynomials of degree 4 with up to six variables.  Again it appears to
be not an  easy task to find all local minima reliably and quickly. Thus
without introduction of new ideas we have the impression that the
algorithm is practically restricted to 3-space.

\section{Necessary conditions for test sets} \label{num_red}

In the following let $P\in {\cal K}_0^3$ be a centrally symmetric polytope with $n$ facets $F_i$ and let $H_i=\aff(F_i)$, $1\leq i \leq n$. By $H^+_i$ we denote the halfspace bounded by $H_i$ containing $P$ and let $H^-_i=\R^3\backslash H^+_i \cup H_i$. We always assume that we have a {\em lattice description}  of the polytope, i.e., we know the face lattice of $P$ specified by its Hasse diagram and  the vertices and edges of $P$ (cf.~\cite{Seidel:computational}). In particular, for each facet $F_i$ we have a list ${\cal N}(F_i)$ of its edge-neighbours, i.e., 
$$
 {\cal N}(F_i) =\{ F_j: \dim(F_i\cap F_j)=1\}.
$$  
We remark that such a lattice description can be computed from the description $P=\cap_{i=1}^n H^+_i$ in time $O(n\log n)$ \cite{Chazalle:convexhull}. Regarding  the combinatorics of polytopes  we refer to the books \cite{McMullenShephard:polytopes} and \cite{Ziegler:polytopes}. 

\vspace{0.1cm}

As pointed out in the last section one crucial point of Algorithm \ref{algorithm_org} is the number of choices of facets (or hyperplanes $H_i$) which have to be considered for each case of the algorithm. In this section we show how one can reduce this number. However, since it turns out that the most time consuming step of Algorithm \ref{algorithm_org} is step S2., we are also looking for ways to reduce the number of executions of step S2.~as well.  An exact worst case analysis of
the complexity of the resulting algorithm appears to be not completely
straightforward. But we show that for  some rather natural
classes of polytopes with $n$ facets we can eliminate ``most'' possible choices of facets in time $O(n^{2})$ 
and we have to carry out only $O(n^{3/2})$ times the steps  S1., S2.~and S3.

Of course, using the central symmetry of the polytope and   
 the arbitrariness of the order of the basis of a lattice we can reduce the number of  possible choices of hyperplanes $H_i$ in the first two cases   to $\frac{1}{12}\binom{n}{6}$   and to 
$2\cdot \frac{1}{12} \binom{n}{7}$ choices in the remaining cases. This does not really help and so one could try to make further use of the symmetries of a given polytope,  
 as {\rm Minkowski} did in his study of the octahedron, where
he managed to reduce the number of cases to 1 (!). Thus he did not need to carry out one step S2. 
However, for polytopes with  little symmetry this would not be of great help  and 
therefore we use a different approach.

With respect to step S3.~of algorithm \ref{algorithm_org} we are only interested in a selection of hyperplanes $H_i$, $1\leq i \leq k$, (say)  such that 
\begin{equation}
\label{feasy}
  {\cal S}_{H_1,\dots,H_k} \cap \{W \in \R^{3\times 3} : {\cal U}_W\subset \bd(P)\}\ne\emptyset,  
\end{equation} 
where ${\cal U}_W$ is a test set of the first, second or third kind corresponding to the case we are studying. If ${\cal U}_{W}=\{u^1_W,\dots,u^k_W\}$ then \eqref{feasy} just says that  the point $u_W^i$ should not only lie in the hyperplane $H_i$, but in the facet $F_i$. Hence \eqref{feasy} can be reformulated as the condition that 
\begin{equation} 
\label{feaslp}
\begin{split}
 {\widetilde {\cal S}}_{H_1,\dots,H_K}  = 
\Big\{ W\in \R^{3\times 3}  & : u^i_W \in H_i \text{ and } 
  u^i_W \in H^+_{i_j}  \\ &\quad\text{for all } F_{i_J}\in{\cal N}(F^i),\;1\leq i \leq k  \Big\} \ne\emptyset. 
\end{split}
\end{equation}  
To check whether ${\widetilde {\cal S}}_{H_1,\dots,H_K}$ is empty is just an instance of a feasibility problem of Linear Programming and may be easily solved by any LP-solver (cf e.g.~\cite{Schrijver:programming}). 

Since the vectors in a test set are linearly dependent the set ${\widetilde {\cal S}}_{H_1,\dots,H_K}$ will be empty for ``most choices'' of hyperplanes $H_i$. For instance: Let $\sigma=-1$ if we are dealing with case 1 and let $\sigma=1$ otherwise.  Then the vectors of a test set satisfy the relations    
\begin{equation}
\label{relations}
u^1_W+\sigma u^2_W=u^6_W,\;\; u^2_W+\sigma u^3_W=u^4_W,\;\; u^1_W+\sigma u^3_W=u^5_W.
\end{equation}
If we have fixed a facet $F_1$ with hyperplane $H_1=\{x\in \R^3 : (a^1)^\intercal x=b_1\}$, say, for the vector $u^1_W$, then by the first relation of \eqref{relations} we get   $\sigma (a^1)^\intercal u^2_W \leq 0$. Otherwise the sum of these two vectors would be separated by $H_1$ from the polytope, but  the sum $u^6_W$ has to lie in a facet. This trivial observation already reduces the possible choices for $H_2$ to almost $n/2$ and obviously we  can apply the same argumentation to the other hyperplanes. 
However, as we shall see at the end of the section, a detailed analysis of the linear dependencies will give a much better reduction for certain classes of polytopes. 

Of course, the determination of the emptiness of \eqref{feaslp} would reduce the number of executions of step S2, but we still have to consider all possible choices. Hence we have the problem  to find a ``fast'' way to exclude  ``almost'' all choices of facets (hyperplanes $H_i$) with  ${\widetilde {\cal S}}_{H_1,\dots,H_K}=\emptyset$. To this end we make  use of \eqref{relations}. Let 
\begin{equation}
\label{set_G}
 {\cal G}  = \left\{ (F_i,F_j,F_k) : (F_i + \sigma F_j) \cap  F_k \ne \emptyset \right\}.
\end{equation}
Obviously, 
\begin{equation*}
{\widetilde {\cal S}}_{H_{l_1},\dots,H_{l_k}}\ne \emptyset \Longrightarrow 
 (F_{l_1},F_{l_2},F_{l_6}),  (F_{l_2},F_{l_3},F_{l_4}),  (F_{l_1},F_{l_3},F_{l_5}) \in {\cal G}, 
\end{equation*}
 and to test whether a tuple $(F_i,F_j,F_k)$ belongs to ${\cal G}$ is just a feasibility problem of Linear Programming, namely:
\begin{equation}
\label{lp_g}
\begin{split}
 (F_i,F_j,F_k)\in {\cal G} \Leftrightarrow  &
\Big \{  (w^1,w^2) \in \R^{3 \times 2}  :  w^1 \in H_i,\; w^1 \in H_{i_l}^+ \text{ for } F_{i_l}\in {\cal N}(F_i), \\ 
& w^2 \in H_j,\; w^2 \in H_{j_l}^+ \text{ for } F_{j_l}\in {\cal N}(F_j), \\
& w^1+\sigma w^2 \in H_k,\; w^1+\sigma w^2 \in H_{k_l}^+ \text{ for } F_{k_l}\in {\cal N}(F_k)
 \Big\} \ne \emptyset. 
\end{split}
\end{equation}

The construction of the set ${\cal G}$ by enumeration clearly involves the consideration of $O(n^3)$ possibilities. While in the numerical examples given in the last section the most time consuming part of the algorithm was the solution of (\ref{feaslp}), the series of examples at the end of the section indicates that  for large polytopes the determination of ${\cal G}$ could be the hardest part. Therefore we show now how the geometry of the polytope can be used to do this effectively. To this end let 
\begin{equation*}
  {\cal G}(F_i)=\left\{ F_j : (F_i+\sigma F_j) \cap \bd(P) \ne\emptyset \right\},\quad 1\leq i\leq n. 
\end{equation*}
With this notation we have 

\begin{lemma}\label{boundary}
Let $F_i$ be a facet of $P$ and let $F_i^\cup=\bigcup_{F_j\in {\cal G}(F_i)}F_j$. The set $F_i^\cup$ is edge-connected, i.e., any two points $x,y\in F_i^\cup$ can be connected by a continuous path contained in $ F_i^\cup$ without crossing a vertex of $P$.
\end{lemma}
\begin{proof}
For a fixed point $v\in F_i$ let ${\cal I}(v)=\bd(P) \cap \sigma (\bd(P) -v)= \bd(P) \cap (\bd(P) -\sigma v)$. Then we have  ${\cal I}(v)=\{y\in P\cap (P-\sigma v): y+\mu v,\, \mu\in\R, \text{ is a supporting line of }P\cap (P-\sigma v)\}$ and  thus ${\cal I}(v)$ is  the shadow boundary of $P\cap (P-\sigma v)$ in direction $v$. Hence  two points of ${\cal I}(v)$ can be connected by a continuous  path. Now let 
$$
  {\cal I}(v)^\cup=\cup_{\{F_j: F_j\cap {\cal I}(v)\ne \emptyset\}} F_j.
$$
Then, by construction, the set ${\cal I}(v)^\cup$ is edge-connected and so it is the union 
$\bigcup_{v\in F_i}{\cal I}(v)^\cup=F_i^\cup$. 
\end{proof}
Once we have found one element of the set ${\cal G}(F_i)$, the lemma says that we  can determine the other elements by recursively checking the neighbours of the elements which were already found and since 
\begin{equation}
\label{do we need it}
(F_i,F_j,F_k) \in {\cal G} \Rightarrow F_j  \in {\cal G}(F_i),  
\end{equation}
we may use Lemma \ref{boundary} in order to construct the set ${\cal G}$. In order 
to present an algorithm for computing ${\cal G}$ as well as the sets ${\cal G}(F_i)$ we need one more notation: For $l,u \in \R^3$, $l \leq u$,  let $B(l,u)=\{x\in \R^3 : l\leq x \leq u\}$ be the box parallel to the coordinate axes with lower vertex $l$ and upper vertex $u$ and for a facet $F_i$ let $B(l^i,u^i)$ be the minimal box containing $F_i$, which can be  computed from the coordinates of the vertices of $F_i$ in time $O(n)$. A necessary condition for $(F_i+\sigma F_j)\cap F_k\ne \emptyset$ is given by  
\begin{equation}
\label{boxsum}
\begin{split}
\big(B(l^i,u^i)+&\sigma B(l^j,u^j)\big) \cap B(l^k,u^k) = \\ 
& \begin{cases}
B(l^i+l^j,u^i+u^j)\cap B(l^k,u^k) \ne \emptyset, &\sigma=1, \\
B(l^i-u^j,u^i-l^j)\cap B(l^k,u^k) \ne \emptyset, & \sigma=-1,
\end{cases}
\end{split}
\end{equation}
which can  easily be checked since  
\begin{equation}
\label{boxintersection}
 B(l,u)\cap B(\tilde{l},\tilde{u})\ne\emptyset \Longleftrightarrow l\leq \tilde {u} \text{ and } \tilde{l}\leq u.
\end{equation}
Therefore, if we want to find for a given pair $F_i,F_j$ all facets $F_k$ with $(F_i,F_j,F_k)\in {\cal G}$ we just have to consider the facets corresponding to boxes $B(l^k,u^k)$ satisfying \eqref{boxsum}. For polytopes with many facets ``most'' facets will be ``far away'' from $F_i+F_j$ and thus \eqref{boxsum} won't be fulfilled for ``most'' facets. 

The next lemma says how we can find  one facet lying in a set ${\cal G}(F_i)$ with the help of these boxes .
\begin{lemma} 
\label{one_facet}
Let $\nu(P)=\max\{\#{\cal N}(F_i): 1\leq i\leq n\}$, let $\beta(P)$ be the maximal number of boxes $B(l^k,u^k)$ intersecting a box of the form $B(l^i+u^i)+\sigma B(l^j+u^j)$, $1\leq i,j \leq n$, and let $\eta(P)$ be the maximal number of facets of $P$ having a nonempty intersection with a fixed hyperplane containing the origin. Then a facet $F_j$ lying in a set ${\cal G}(F_i)$ can be found in time $O(n+\eta(P)\log^2(n)+\eta(P)\beta(P)\nu(P))$.
\end{lemma}
\begin{proof}
Let $v\in F_i$ and let $H$ be any hyperplane containing $v$ and the origin. Then there exists a facet $F_j\in{\cal G}(F_i)$ having a nonempty intersection with $H\cap P$. Using the edges of $P$ we can easily determine all facets $\{F_{j_1},\dots, F_{j_l}\}$, $j_l\leq \eta(P)$, with this property. Since the polytope has  $O(n)$ edges this can be carried out in time $O(n)$. For each $j_i$ the boxes $B(l^k,u^k)$, $k\leq \beta(P)$,   intersecting $\left(B(l^i,u^i)+\sigma B(l^{j_i},u^{j_i})\right)$ can be determined in time $O(\log^2(n)+\beta(P))$ 
by  well know methods from computational geometry about range searching (cf.~\cite{Agarwal:handbook}). For each possible choice $F_i, F_{j_i}, F_k$ we use {\rm \eqref{lp_g}} to verify whether  $(F_i,F_{j_i}, F_{k})\in {\cal G}$ and thus $F_{j_i}\in{\cal G}(F_i)$. Now each \eqref{lp_g} is a feasibility problem of Linear Programming  with $O(\nu(P))$ constraints in dimension $6$ and this can be solved  with $O(\nu(P))$ arithmetic operations (cf.~\cite{Schrijver:programming}, pp.~199).
\end{proof}

Using the last lemma we have the following algorithm for computing the set ${\cal G}$ and the sets ${\cal G}(F_i)$, $1\leq i \leq n$.

\begin{algorithm}[Computing ${\cal G}(F_i)$ and ${\cal G}$] \hfill
\label{algorithm_g}
\newline {\rm Input:} A polytope $P\in {\cal K}_0^3$ given by the supporting hyperplanes $H_i$, $1\leq i\leq n$, of its facets $F_i$, as well as a lattice description of $P$. 
\hfill\newline
{\rm Output:} ${\cal G}$ and ${\cal G}(F_i)$, $1\leq i\leq n$. 
\vspace{0.0cm}
\begin{itemize}
\item  Let $U=\emptyset$ and determine the boxes $B(l^i,u^i)$, $1\leq i \leq n$.
\item For each facet $F_i$ do \vspace{-0.2cm}
\begin{enumerate} 
\item[{\rm G1.}] Find a facet $F_j\in {\cal G}(F_i)$. 
                Let $U_i=\emptyset$ and $N=\{F_j\}$. 
\item[{\rm G2.}] While $N\ne\emptyset$ do 
    \addtolength\itemindent{0.2cm}
\item[{\rm G3.}] For $F_{j_1}\in N$ determine all facets $F_{k_1},\dots, F_{k_p}$ such                 that \hfill\newline $\left(B(l^i,u^i)+\sigma B(l^{j_1},u^{j_1})\right)\cap B(l^{k_l},u^{k_l}) \ne \emptyset$ {\rm (cf.~\eqref{boxsum}, \eqref{boxintersection})}.
\item[{\rm G4.}] Use {\rm \eqref{lp_g}} to determine all $(F_i,F_{j_1}, F_{k_l})\in {\cal G}$, 
                $k_l\in\{k_1,\dots,k_p\}$, and add them to the set $U$.
\item[{\rm G5.}] Let $U_i=U_i\cup\{F_{j_1}\}$ and $N=(N\cup{\cal N}(F_{j_1}))\backslash U_i$, if there exists a $(F_i,F_{j_1}, F_{k_l})\in {\cal G}$, otherwise let $N=N\backslash F_{j_1}$. 
  \end{enumerate}
\end{itemize}
\end{algorithm}

\begin{lemma}  
\label{prep_lemma}
Let $\gamma(P)=\max\{\#{\cal G}(F_i): 1\leq i\leq n\}$,  and let $\nu(P)$, $\eta(P)$, $\beta(P)$  as in Lemma {\rm \ref{one_facet}}.  Then Algorithm \ref{algorithm_g} determines ${\cal G}$ and ${\cal G}(F_i)$, $1\leq i\leq n$, in time 
$O\big(n^2+n\log^2(n)(\eta(P)+\gamma(P)\nu(P)\beta(P))+n\beta(P)\nu(P)(\eta(P)+\allowbreak \gamma(P)\nu(P))\big)$.
\end{lemma}

\begin{proof}
By Lemma \ref{boundary} the algorithm finds all $(F_i,F_j,F_k)\in {\cal G}$ and at the end we have $U={\cal G}$. Furthermore,  at the end of each loop G2.~the set $U_i$ coincides with ${\cal G}(F_i)$. 

To find a first facet $F_j\in{\cal G}(F_i)$ we need by Lemma \ref{one_facet} at most $O(n+\eta(P)\log^2(n)+\eta(P)\beta(P)\nu(P))$ operations. For the estimation of the steps G3.~and G4.~we can proceed as in the proof of Lemma \ref{one_facet}. The boxes in step G3.~can be found  in time $O(\log^2(n)+\beta(P))$ and in step G4.~ we have to solve at most $\beta(P)$ feasibility problems (cf.~\eqref{lp_g}) with at most $\nu(P)$ constraints in dimension 6 which can be done with $O(\nu(P))$ arithmetic operations. Finally, we observe that for each facet $F_i$ the  loop G2.~is executed at most $\gamma(P)\nu(P)$ times. 
\end{proof}
\par\noindent\textbf{Remark.} 
\begin{enumerate}
\item[{\rm i)}] The bound on the running time given in the last lemma is useless for a worst-case analysis, because there exist polytopes such that each of the numbers $\eta(P),\beta(P),\gamma(P),\nu(P)$ is of order $O(n)$. In this case Lemma \ref{prep_lemma} would give a $O(n^5)$-algorithm and of course, one can determine the sets ${\cal G}, {\cal G}(F_i)$ by a trivial $O(n^4)$-algorithm. 
 Nevertheless we shall see in Theorem \ref{theorem_one} that  this lemma gives a $O(n^2)$ bound for a rather natural class of polytopes. 
\item[{\rm ii)}] As we have $F_j\in{\cal G}(F_i)\Leftrightarrow F_i\in {\cal G}(F_j)$ we may use the  facets belonging to a set ${\cal G}(F_i)$ as starting facets in step G1. It is not hard to see that we can determine all starting facets except the first one in this way. However, for simplification we do not exploit this fact. 
\end{enumerate}

\vspace{0.3cm}
\par\noindent
Altogether the previous observations lead to the following refinement of Algorithm \ref{algorithm_org}
\begin{algorithm} \hfill
\label{algorithm_2}
\par\noindent {\rm Input:} A polytope $P\in {\cal K}_0^3$ given by the supporting hyperplanes $H_i$, $1\leq i\leq n$, of its facets $F_i$, as well as a lattice description of $P$.  \hfill 
\par\noindent {\rm Output:} A densest packing lattice.
\vspace{-0.2cm}
\begin{itemize}
\item For each of the  cases {\rm I., II., III., IV.}~of Theorem {\rm \ref{ness_con}} do \vspace{-0.2cm}
\begin{itemize}
\item [{\rm I1.}] Compute the sets ${\cal G}(F_i)$ and  ${\cal G}$ with Algorithm {\rm \ref{algorithm_g}}.
\item For three facets  $F_{l_1},F_{l_2},F_{l_3}$   satisfying $F_{l_2} \in{\cal G}(F_{l_1})$ and $F_{l_3}\in {\cal G}(F_{l_1})\cap {\cal G}(F_{l_2})$ do  
    \addtolength\itemindent{0.2cm}

\item For every choice of  facets $F_{l_i}$, $4\leq i\leq k$,  with $(F_{l_1},F_{l_2},F_{l_6}), \allowbreak  (F_{l_2},F_{l_3},F_{l_4}),  (F_{l_1},F_{l_3},F_{l_5}) \in {\cal G}$ do 
\begin{enumerate}
    \item[{\rm S0.}] If ${\widetilde {\cal S}}_{H_{l_1},\dots,H_{l_k}}\ne \emptyset$ {\rm (cf.~\eqref{feaslp})} do 
    \addtolength\itemindent{0.2cm}
    \item[{\rm S1.}] Determine ${\cal S}_{H_{l_1},\dots, H_{l_k}}$.
    \item[{\rm S2.}] Find the local minima ${\cal M}_{H_{l_1},\dots,H_{l_k}}$ of $f_{H_{l_1},\dots,H_{l_k}}$. 
    \item[{\rm S3.}] For each  $M\in {\cal M}_{H_{l_1},\dots,H_{l_k}}$ check whether $M\cdot \Z^d$ is an admissible lattice, i.e., $M$ satisfies the criterion {\rm 1}.~of Lemma {\rm \ref{lattice_check}} in the first case and criterion {\rm 2.}~in the remaining cases. 
\end{enumerate}
\end{itemize}
\item Among all calculated admissible lattices find one with  minimal determinant. The corresponding lattice is a critical lattice of $P$.
\end{itemize}
\end{algorithm}
\noindent 
It seems to be a  non trivial problem to give a ``non trivial''  worst case analysis of the algorithm for an arbitrary polytope. In the following we want to demonstrate, by a rather natural series of polytopes,  the improvement of Algroithm \eqref{algorithm_2} compared with  the brute force method (see Algorithm \eqref{algorithm_org}) which involves the examination of $\Omega(n^7)$ steps S0., S1., S2., and S3.

For a facet $F$ of a polytope $P$ we denote by $R(F)$ its circumradius and by $r(F)$ its inradius with respect to its affine hull. For a real $c>0$ we say  that the facets of $P$ are of $c$--uniform shape, if 
$$
\min\left\{ r(F) : F \text{ facet of }P\right\} \ge c \cdot
\max\left\{ R(F) : F \text{ facet of }P\right\}.
$$

\begin{theorem}\label{theorem_one}
Let $\{P_m\}_{m\in\N}$ be a series of polytopes $P_m\in{\cal K}_0^3$ such that {\rm i)} all $P_m$ have facets of $c$--uniform shape for some fixed $c$ and {\rm ii)} $\{P^m\}$ converges to the 3-dimensional unit ball with respect to the Hausdorff metric. Let $f_m$ be the number of facets of the polytope $P_m$. Then the number  of all choices of facets which will be examined by Algorithm \eqref{algorithm_2} in the steps S0.,S1., S2.~and S3.~ is of size $O(f_m^{3/2})$ and these possible choices of facets can be determined in time $O(f_m^{2})$.
\end{theorem}
\begin{proof}
We project  $\bd(P_m)$ by the radial projection $\rho$  onto the unit sphere $S^2$ and we carry out  all calculations on $S^2$. It can easily be checked that our asymptotical estimates remain correct. The facets of $P_m$ will be denoted by $F^m_i$, $1\leq i \leq f_m$, and in the following we shall denote by $c_i$ certain positive constants.  First we observe that the uniformity of our sequence implies that  for the spherical diameter $d(\rho(F^m_i))$ and the spherical area $A(\rho(F_i^m))$  of the facets holds  
\begin{equation} \label{spher_diam}
c_1/\sqrt{f_m} \le d(\rho(F^m_i)) \le c_2/\sqrt{f_m},\quad  c_3/f_m \le A(\rho(F^m_i)) \le c_4/f_m.
\end{equation}
Next we note that for a point $x \in S^2$ we have
\begin{equation}\label{belt}
\{y \in S^2 : y+x \in S^2 \} = \{ y \in S^2 : \langle x,y \rangle = -1/2\}.
\end{equation}
Now let $F^m_i$, $F^m_j$ be two facets of $P_m$ such that $(F^m_i+\sigma F^m_j)\cap \bd(P_m) \ne \emptyset$ and let $x\in F_i^m$. We find by (\ref{spher_diam}), (\ref{belt}) that  
\begin{equation} \label{big_belt}
F^m_j \subset \{ y \in  \bd(P_m) : -\sigma/2-c_5/\sqrt{f_m} \le \langle \rho(y),\rho(x) \rangle \le - \sigma/2+c_5/\sqrt{f_m}\}.
\end{equation}
 As the spherical area  of the set on the right is bounded by $c_6/\sqrt{f_m}$ we find by (\ref{spher_diam}) that for fixed $F^m_i$ there are at most $O(\sqrt{f_m})$ faces $F^m_j$ such that $(F^m_i+\sigma F^m_j) \cap \bd(P_m) \neq \emptyset$. For two faces $F^m_i$ and $F^m_j$ we have $R(F^m_i+\sigma F^m_j) \le R(F^m_i)+R(F^m_j) \le c_7/\sqrt{f_m}$. Proceeding as before we see  that there can be at most $O(1)$ faces of $P_m$ which intersect $(F^m_i+\sigma F^m_j)$. Altogether we have found that 
\begin{equation}
\label{cardinal}
 \#{\cal G}=O(f_m^{3/2})\quad \text{and}\quad \#{\cal G}(F^m_i)=O(\sqrt{f_m}).
\end{equation}
Next we  ask for  the efficiency of Algorithm \ref{algorithm_g}  to determine these sets.  Let $B(l^{m,i},u^{m,i})$ be a minimal box containing $F^m_i$. By the above estimate for the circumradius of two facets $F^m_i,F^m_j$ we have  $|l^{m,i}_k+l^{m,j}_k-u^{m,i}_k-u^{m,j}_k|=O(1/\sqrt{f_m})$ for each coordinate $k$. Again using the area bound \eqref{spher_diam} we see that there are at most $O(1)$ facets $F^m_s$ with $u^{m,s}_1\geq l^{m,i}_1+l^{m,j}_1$ (or $u^{m,s}_1\geq l^{m,i}_1-u^{m,j}_1$) and $l^{m,s}_1\leq u^{m,i}_1+u^{m,j}_1$ (or $l^{m,s}_1\leq u^{m,i}_1-l^{m,j}_1$ ) (cf.~\eqref{boxintersection}). Hence the number $\beta(P)$ of Lemma \ref{prep_lemma} is of order $O(1)$. Next we note that on account of our assumptions the maximal number $\nu(P)$ of  neighbours  of a given facet is constant, if $f_m$ is large enough. Moreover, by \eqref{big_belt} we also see that the number $\eta(P)$ is of order $\sqrt{f_m}$ and 
therefore Algorithm \ref{algorithm_g} determines the set ${\cal G}$ as well as the sets ${\cal G}(F^m_i)$ in time $O(f_m^{2})$ (cf.~Lemma \eqref{prep_lemma}).

Since we have already proven that $\#{\cal G}(F^m_i)=O(\sqrt{f_m})$ and that for two given facets $F^m_{l_1}, F^m_{l_2}$ there are only  $O(1)$ many facets $F^m_{l_k}$ with $(F^m_{l_1}, F^m_{l_2},F^m_{l_k})\allowbreak\in{\cal G}$, it remains to show that   $\#({\cal G}(F^m_{l_1})\cap{\cal G}(F^m_{l_2}))=O(1)$ for $F^m_{l_2}\in{\cal G}(F^m_{l_1})$ (cf.~Algorithm \ref{algorithm_2}).
Let $x^i \in {\cal G}(F^m_{l_i})$. Then we have (cf.~\eqref{big_belt})
\begin{equation*}
\begin{split}
& -\frac{\sigma}{2}-\frac{c_5}{\sqrt{f_m}} \le \langle \rho(x^1),\rho(x^2) \rangle \le -\frac{\sigma}{2}+\frac{c_5}{\sqrt{f_m}}, \text{ and }\quad {\cal G}(F^m_{l_1})\cap{\cal G}(F^m_{l_2})  \\
& \subset \left\{ y \in  \bd(P_m) : - \frac{\sigma}{2}-\frac{c_5}{\sqrt{f_m}} \le \langle \rho(y),\rho(x^i) \rangle \le - \frac{\sigma}{2}+\frac{c_5}{\sqrt{f_m}}, i=1,2\right\}.
\end{split}
\end{equation*}
The radial projection of the latter set is a spherical parallelogram with ``edge length'' $O(1/\sqrt{f_m})$. Hence its spherical area is $O(1/f_m)$ and together with \eqref{spher_diam} this shows $\#({\cal G}(F^m_{l_1})\cap{\cal G}(F^m_{l_2}))=O(1)$. 

\end{proof}

\par\noindent \textbf{Remark.}~Obviously, in the cases II, III.~and  IV.~the sets  ${\cal G}(F_i)$, ${\cal G}$ coincide. Furthermore, since $(F_i-F_j)\cap F_k =(F_i+(-F_j))\cap F_k$ and $-F_j$ is a facet of $P$ it suffices to determine ${\cal G}(F_i)$, ${\cal G}$ only for one case.

\section{Determination of the local minima} \label{min_pol}

In this section we concentrate on the steps S1., S2.~and S3.~of Algorithm \ref{algorithm_2}. To this end let $P\in {\cal K}_0^3$ be a centrally symmetric polytope with $n$ facets given by the inequalities 
\begin{equation}
 P=\left\{ x \in \R^3 : (a^i)^\intercal x \leq b_i,\, 1\leq i\leq n\right\}, 
\end{equation}
where $a^i\in \R^3$ and $b_i\in \R_{>0}$. Let $H_i=\{x \in \R^3 : (a^i)^\intercal x = b_i\}$ be the supporting hyperplane of the facet $F_i$. Since the cases I -- IV of the algorithm can be treated `` more or less'' in the same way, in the following we shall focus only on case {\rm I} and at the end of the section we shall discuss the necessary changes and adoptions for the other cases. 

For six hyperplanes $H_{i_j}$, $1\leq j\leq 6$, the set ${\cal S}_{H_{i_1},\dots,H_{i_6}}$ (see step S1.~and \eqref{set_S}) is given by all $W=(w^1,w^2,w^3)\in \R^{3\times 3}$ satisfying 
\begin{equation}
\label{matrix}
\begin{pmatrix}
  (a^{i_1})^\intercal &                 0 &              0     \\
  0                 & (a^{i_2})^\intercal &              0     \\
  0                 &                 0 & (a^{i_3})^\intercal  \\ 
  0                 & (a^{i_4})^\intercal & -(a^{i_4})^\intercal  \\
 -(a^{i_5})^\intercal &                 0 & (a^{i_5})^\intercal   \\
  (a^{i_6})^\intercal & -(a^{i_6})^\intercal &              0     
\end{pmatrix}  
\, 
\begin{pmatrix}
 w^1 \\ w^2 \\ w^3 
\end{pmatrix}
\,
=
\,
\begin{pmatrix}
   b_{i_1} \\    b_{i_2} \\     b_{i_3} \\  
   b_{i_4} \\     b_{i_5} \\     b_{i_6} 
\end{pmatrix}.
\end{equation}
We denote the matrix on the left $A_{i_1,\dots,i_6} \in \R^{6\times 9}$ and with suitable matrices $C_{i_1,\dots,i_6}$, $M^j_{i_1,\dots,i_6}$, $1\leq j \leq 9-{\rm rank}(A_{i_1,\dots,i_6})$ we may write  
\begin{equation}
\label{para_set}
  {\cal S}_{H_{i_1},\dots,H_{i_6}} =\left\{ W\in \R^{3\times 3} : 
    W=  C_{i_1,\dots,i_6} +\hspace{-0.9cm} \sum_{j=1}^{9-{\rm rank}(A_{i_1,\dots,i_6})} 
    \lambda_j \cdot M^j_{i_1,\dots,i_6},\,\, \lambda_j \in \R 
   \right\}. 
\end{equation}
We remark that the matrices $C_{i_1,\dots,i_6}$, $M^j_{i_1,\dots,i_6}$ can easily be determined by any program for solving systems of linear equations.


\begin{lemma} 
\label{rang_lemma}
If ${\rm rank}(A_{i_1,\dots,i_6})<6$ then  the function $f_{H_{i_1},\dots, H_{i_6}}(W)$ {\rm (}cf. \eqref{function}{\rm )}  has no local  minimum $W\in {\cal S}_{H_{i_1},\dots,H_{i_6}}$ with $ f_{H_{i_1},\dots, H_{i_6}}(W)>0$.
\end{lemma} 
\begin{proof} Since the vectors $a^{j_i}$ correspond to facets defining hyperplanes the vectors $(a^{i_1},0,0)^\intercal, (0,a^{i_2},0)^\intercal, (a^{i_6},-a^{i_6},0)^\intercal \in \R^9$ are linearly independent. 
Otherwise we can assume that $a^{i_6}=a^{i_1}=\pm a^{i_2}$ and we get $(a^{i_6})^\intercal (w^1-w^2)=b_{i_1}\mp b_{i_2}\ne b_{i_1}=b_{i_6}$. 
Hence the vectors $\{a^{i_3}, a^{i_4}, a^{i_5}\}$ are linearly dependent, because otherwise ${\rm rank}(A_{i_1,\dots,i_6})=6$.  In the same way we find that $\{a^{i_1}, a^{i_5}, a^{i_6}\}$ and $\{a^{i_2}, a^{i_4},\allowbreak a^{i_6}\}$ are linearly dependent. Thus we can find three nontrivial vectors $v^1,\allowbreak v^2,\allowbreak v^3 \in \R^3\backslash\{0\}$ such that 
$$
  v^1\in \lin \{a^{i_3}, a^{i_4}, a^{i_5}\}^\perp, \, 
  v^2\in \lin \{a^{i_1}, a^{i_5}, a^{i_6}\}^\perp, \,
  v^3\in \lin \{a^{i_2}, a^{i_4}, a^{i_6}\}^\perp,
$$
where $\lin U^\perp$ denotes the orthogonal complement of $\lin U$. 
Now let $W=(w^1,w^2,w^3) \in {\cal S}_{H_{i_1},\dots,H_{i_6}}$ with $\det(W)\ne 0$, and let us assume that $W$ is a local minimum.  By the choice of $v^i$ we have $(w^1+\lambda v^1,w^2+\mu v^2, w^3+\nu v^3)\in     {\cal S}_{H_{i_1},\dots,H_{i_6}}$ for all $\lambda,\mu, \nu\in\R$ and therefore let 
\begin{align*}
 g(\lambda,\mu,\nu)&=\det(w^1+\lambda v^1,w^2+\mu v^2, w^3+\nu v^3)=\det(W) \\
     &+\, \lambda\, \det(v^1,w^2,w^3) + \mu\, \det(w^1,v^2,w^3) + \nu\, \det(w^1,w^2,v^3)\\
     &+\, \lambda\mu\,\det(v^1,v^2,w^3) + \mu\nu\,\det(w^1,v^2,v^3) + \lambda\nu\, \det(v^1,w^2,v^3) \\
&+\, \lambda\mu\nu \det(v^1,v^2,v^3). 
\end{align*}
It is easy to see that this function has a local extremum at $(0,0,0)$ if and only if it is constant, i.e., $g(\lambda,\mu,\nu)=\det(W)$. In particular  we have $\det(v^1,w^2,w^3)=\det(w^1,v^2,w^3)=\det(v^1,v^2,w^3)=0$. 
Since $w^1,w^2,w^3$ are linearly independent this implies that $v^1$ or $v^2$ belongs to $\lin\{w^3\}$, and in the same way we find that $v^2$ or $v^3$ lies in  $\lin\{w^1\}$ and $v^1$ or $v^3$ belongs to $\lin\{w^2\}$. However, since $v^i\in \R^3\backslash\{0\}$ this yields the contradiction $|\det(v^1,v^2,v^3)|=|\det(W)|$.

\end{proof}
Since we are only interested in local minima $W$ of the functions $f_{H_{i_1},\dots, H_{i_6}}$ with $f_{H_{i_1},\dots, H_{i_6}}(W)> 0$, in the following we assume that ${\rm rank}(A_{i_1,\dots,i_6})=6$. Instead of searching for the local minima of  $f_{H_{i_1},\dots, H_{i_6}}(W)=|\det(W)|$  it is more practical to look for all local extrema of the function $\det(W)$, which can be parameterized by ${\mathfrak p}_{i_1,\dots,i_6}: \R^3 \to \R$ given by  (see \eqref{para_set})
\begin{equation}
\label{our_polynomial}
  {\mathfrak p}_{i_1,\dots,i_6}(x,y,z) =
  \det\left( C_{i_1,\dots,i_6} +x\cdot M^1_{i_1,\dots,i_6}+
             y\cdot M^2_{i_1,\dots,i_6}+z\cdot M^3_{i_1,\dots,i_6} \right).
\end{equation}
$\mathfrak{p}_{i_1,\dots,i_6}$ is a polynomial in $3$ variables $(x,y,z)$, where each monomial has total degree $3$  at most. In what follows we are mainly interested in general properties of such a polynomial, and therefore we shall write $\mathfrak p$ instead of ${\mathfrak p}_{i_1,\dots,i_6}$. So $\mathfrak{p}$ can be written as 
$$
\mathfrak{p}(x,y,z)= \sum_{0\leq i,j,k\leq 3, \,\,\,  i+j+k\leq 3} \alpha_{i,j,k}\cdot x^iy^jz^k, 
$$
for some scalars $\alpha_{i,j,k}\in \R$. The canonical first step in order to find the local extrema is to calculate the set ${\cal V}(\nabla\mathfrak p) $ where the gradient $\nabla\mathfrak{p}$ vanishes, i.e.,
\begin{equation*} 
 {\cal V}(\nabla\mathfrak p)=\left\{ (x,y,z)\in \R^3 : \nabla\mathfrak{p}(x,y,z)=0 \right\}.
\end{equation*}
Thus we are interested in the common roots of the partial derivatives
\begin{equation}
\label{partial}
\begin{split}
 \mathfrak{p}_1 & =  \partial\mathfrak{p}/\partial x= \chi_1 x^2 + \mathfrak{l}_1(y,z)\cdot x  + 
                  \mathfrak{q}_1(y,z), \\
 \mathfrak{p}_2 & =  \partial\mathfrak{p}/\partial y= \chi_2 x^2 + \mathfrak{l}_2(y,z)\cdot x + 
                  \mathfrak{q}_2(y,z), \\
 \mathfrak{p}_3 & =  \partial\mathfrak{p}/\partial z= \chi_3 x^2 + \mathfrak{l}_3(y,z)\cdot x + 
                  \mathfrak{q}_3(y,z),
\end{split}
\end{equation} 
where $\chi_i\in \R$, $\mathfrak{l}_i=\chi_{i,2}y + \chi_{i,3}z+\chi_{i,0}$, $\chi_{i,j}\in \R$, and $\mathfrak{q}_i=\chi_{i,2,0}y^2+ \chi_{i,1,1} yz +\chi_{i,0,2}z^2 + \chi_{i,1,0}y+\chi_{i,0,1}z +\chi_{i,0,0}$, $\chi_{i,j,k}\in \R$, $1\leq i\leq 3$. As $\mathfrak{p}$ is a polynomial of total degree of at most 3, ${\cal V}(\nabla\mathfrak p)$ has the following nice property.
\begin{lemma}
\label{only_affine}
Let ${ m}\in{\cal V}(\nabla\mathfrak p)$ be a local extremum of the function $\mathfrak{p}$ and let ${\cal C}\subset {\cal V}(\nabla\mathfrak p)$ be a (path-) connected component containing ${ m}$. Then $\aff({\cal C})\subset {\cal V}(\nabla\mathfrak p)$. 
\end{lemma}
\begin{proof} Let ${ n}\in {\cal C}$, ${ m}\ne{ n}$.  Since there exists a path from ${ m}$ to ${ n}$ in ${\cal C}$ and since the gradient vanishes on ${\cal C}$ we have $\mathfrak{p}({ m})=\mathfrak{p}({ n})$. For $t\in\R$ let $s(t)={ m}+t\cdot ({ n}-{ m})$. The function $\mathfrak{p}(s(t))$ is an univariate polynomial in $t$ of degree at most three with $\mathfrak{p}(s(0))=\mathfrak{p}(s(1))$ and the derivative vanishes at $0$ and $1$. Thus $\mathfrak{p}(s(t))=\mathfrak{p}({ m})$ for all $t\in\R$. Since ${ m}$ is assumed to be a local extremum  there exists a ${\bar t}\in (0,1)$ such that $s({\bar t})$ is a local extremum, too. Hence we have found three points on the line $s(t)$ where the gradient of $\mathfrak{p}$ vanishes. Since all partial derivatives are polynomials of total degree at most two we have shown that $\nabla \mathfrak{p}(s(t)) =0$ for all $t\in \R$.  
\end{proof}

The last lemma tells us that in order to locate all possible local extrema of $\mathfrak{p}$ it suffices to find all {\em isolated  affine subspaces} of $\nabla\mathfrak{p}=0$, i.e., 
\begin{equation}
\label{isolated_subspaces}
\begin{split}
 {\cal V}_\aff(\nabla\mathfrak p)=\Big\{ {\cal C}\subset {\cal V}(\nabla\mathfrak p) : & \,\,   
    {\cal C}=\aff({\cal C}) \text{ and  there exists no connected } \\ 
  & \text{ component } {\cal U}\subset {\cal V}(\nabla\mathfrak p)  \text{ with } {\cal C}  \subsetneqq {\cal U} 
     \Big\}.
\end{split}
\end{equation}  

To our surprise we have not found any efficient algorithm for the determination of the set ${\cal V}$ or ${\cal V}_\aff$ or, in general, for the determination of the local extrema of the polynomial $\mathfrak{p}$. 
 Therefore we have developed an ad hoc method for our purposes which is based on resultants of polynomials and Lemma \ref{only_affine}. 
For a detailed treatment of resultants we refer to \cite{CoxLittleOShea:ideal} and \cite{CoxLittleOShea:using}. Here we just collect some of their basic properties. 

To this end we denote for a polynomial  $\mathfrak{f}\in \R[x_1,\dots,x_d]$ by $\deg(\mathfrak{f})$ its total degree and by $\deg(\mathfrak{f},x_i)$ we denote the degree of the variable $x_i$ if we consider $\mathfrak{f}$ as a polynomial over $\R[x_1,\dots,x_{i-1},x_{i+1},\dots,x_d]$. 
For a polynomial or a system of polynomials ${\cal I} \subset \R[x_1,\dots,x_d]$ we denote by ${\cal V}({\cal I})=\{ x \in \R^n : \mathfrak{f}(x)=0 \text{ for all } \mathfrak{f} \in \cal{I}\}$. Furthermore, let ${\cal V}_\aff({\cal I})$ be the set of all isolated affine subspaces of ${\cal V}({\cal I})$ in the sense of \eqref{isolated_subspaces} and let   
\begin{equation*}
\begin{split}
{\cal V}^j_\aff({\cal I})  = \left\{ {\cal C}\in {\cal V}_\aff({\cal I}) : \dim({\cal C})=j \right\}, \quad j=0,\dots,d. 
\end{split}
\end{equation*}
The elements of ${\cal V}^0_\aff$ will be called {\em isolated roots}.
Since in general is seems to be a hard problem to determine exactly the set  ${\cal V}_\aff({\cal I})$   we only look for  an approximation of this set. This means we want to determine a set 
\begin{equation}
\label{approxaff}
\begin{split}
&\widetilde{\cal V}_\aff({\cal I})\,\, \text{\it consisting of {finitely} many affine subspace }  {\cal C}\subset {\cal V}({\cal I})\\ &\text{\it and }  {\cal V}_\aff({\cal I})\subset \widetilde{\cal V}_\aff({\cal I}).
\end{split}
\end{equation}

Now let $\mathfrak{f}^j\in \R[x_1,\dots,x_d]$, $j=1,2$, be two polynomials and with  respect to the coefficient ring $\R[x_2,\dots,x_d]$ we write 
$$
 f^j(x_1,\dots,x_d) = \sum_{i=0}^{m_j} f_{i,j}\,x_1^i, 
$$
with $f_{i,j}\in \R[x_2,\dots,x_d]$,  $m_j=\deg(\mathfrak{f}^j,x_1)$ and let $m_1, m_2 \geq 1$. 
The {\em resultant} of $f$ and $g$ w.r.t.~$x_1$ will be denoted by $\res(\mathfrak{f}^1,\mathfrak{f}^2,x_1) \in \R[x_2,\dots,x_d]$ and it is given by the determinant of the Sylvester matrix of $\mathfrak{f}^1$ and $\mathfrak{f}^2$ w.r.t.~$x_1$, i.e., 
\begin{equation}
\label{resultant_def}
\res(\mathfrak{f}^1,\mathfrak{f}^2,x_1)=\det
\begin{pmatrix}
    f_{m_1,1}   &         &              & f_{m_2,2}      &               &           \\
    f_{m_1-1,1} &  \ddots &              & f_{m_2-1,2}    &  \ddots       &           \\
    \vdots      &         & f_{m_1,1}    &  \vdots        &               & f_{m_2,2}          \\
    \vdots      &         & f_{m_1-1,1}  &  \vdots        &               & f_{m_2-1,2} \\
    f_{0,1}     &         & \vdots       &                &               & \vdots \\
                &         & \vdots       &   f_{0,2}      &               & \vdots  \\
                &         & f_{0,1}      &                &               &  f_{0,2}
\end{pmatrix},
\end{equation}
where the columns corresponding to $\mathfrak f^1$, $\mathfrak f^2$  are repeated $m_2$-times, $m_1$-times, respectively. The definition can be extended to the case $m_1+m_2\geq 1$ by setting $\res(\mathfrak f^1,\mathfrak f^2,x_1)=\mathfrak f^j$ if $m_j=0$.
\begin{lemma}[{\rm see \cite{CoxLittleOShea:ideal}, pp.~150}] 
\label{resultants_basic}
\par\noindent
\begin{enumerate}
\item[{\rm i)}] Let $(u_1,u_2,\dots,u_d)\in \R^d$ be a common root  of $\mathfrak f^1$ and $\mathfrak f^2$. Then $\res(\mathfrak f^1,\allowbreak\mathfrak f^2, x_1)(u_2,\dots,u_d)=0$.
\item[{\rm ii)}] Let ${ u}=(u_2,\dots,u_d) \in \R^{d-1}$ such that $\res(\mathfrak f^1,\mathfrak f^2,x_1)({ u})=0$ but \linebreak  
      $f_{m_1,1}({ u})\ne 0$ or $f_{m_2,2}({ u})\ne 0$. Then there exists a ${ u_1}\in{ \C}$ such that $\mathfrak f^1({ u_1},{ u})=\mathfrak f^2({ u_1},{ u}) =0$.
\item[{\rm iii)}] $\res(\mathfrak f^1,\mathfrak f^2,x_1)=0$ if and only if $\mathfrak f^1$ and $\mathfrak f^2$ have a common factor $g\in \R[x_1,\dots,x_d]$ with $\deg(\mathfrak g,x_1)\geq 1$. 
\end{enumerate}
\end{lemma}

In general we would like to use resultants in the following way: Let $\mathfrak{p_1}$, $\mathfrak{p_2}$ and $\mathfrak{p_3}$ be the partial derivatives of the polynomial $\mathfrak{p}$ {\rm (see  \eqref{partial})}. First we compute the two resultants $\res_{1,2}=\res(\mathfrak{p}_1,\mathfrak{p}_2,x)\in\R[y,z]$, $\res_{2,3}=\res(\mathfrak{p}_2,\mathfrak{p}_3,x)\in\R[y,z]$ and then the resultant $\res_{1,2,3}=\res\left(\res_{1,2},\res_{2,3},y\right)\allowbreak\in\R[z]$. Next we determine the real roots of $\res_{1,2,3}$, and for each root ${\bar z}$ we determine the common roots of $\res_{1,2}(y,{\bar z})= 0$ and 
$\res_{2,3}(y,{\bar z})=0$. So we get a couple of common roots of $\res_{1,2}$ and $\res_{2,3}$. Again we put each pair of those common roots into $\mathfrak{p}_1,\mathfrak{p}_2,\mathfrak{p}_3$ and solve the three polynomials w.r.t.~$x$.

Now, as we shall see, if $\res_{1,2,3}\ne 0$ (and thus all resultants are non trivial), then we can compute ${\cal V}(\nabla\mathfrak p)$ by this method. However, in general some of the resultants may vanish and hence we have to find common factors of the polynomials (cf.~Lemma \ref{resultants_basic} iii)). According to Lemma \ref{only_affine} we are mainly interested in factors corresponding to affine subspaces. Therefore we call a polynomial $\mathfrak{l}\in \R[x_1,\dots,x_d]$ {\em linear} if $\deg(\mathfrak{l})=1$, i.e.,  
\begin{equation*}
  \mathfrak{l}= l_0+\sum_{i=1}^d l_i \cdot x_i, \quad l_i\in \R. 
\end{equation*}
\begin{remark}
\label{linear_factors}
Let $\mathfrak f=\sum_{i=0}^m f_i\,x_1^i, \mathfrak{l}= l_0+\sum_{i=1}^d l_i \cdot x_i \in R[x_1,\dots,x_d]$ with $f_i\in\R[x_2,\dots,x_d]$, $l_i\in \R$ and $f_m, l_1 \ne 0$. Then $\mathfrak{l}$ is a factor of $\mathfrak f$ if and only if 
$ 
     \mathfrak f\big( ( l_1 x_1- \mathfrak l )/l_1,x_2,\dots,x_d \big) =0. 
$
\end{remark}
\begin{proof} Obviously, if $\mathfrak{l}$ is factor of $\mathfrak f$ then  the statement holds. W.l.o.g.~let $l_1=1$ and $\bar{\mathfrak{l}}=\mathfrak l -  x_1$ and let $\mathfrak{q}=\sum_{i=0}^{m_1} q_i\,x_1^i$ be the polynomial whose coefficients are recursively defined by 
\begin{equation}
\label{rest}
    q_{m-1}= f_m \text{ and } q_i=f_{i+1}-q_{i+1}\bar\mathfrak{l},\quad i=m-2.\dots,0.
\end{equation}
Multiplication of $\mathfrak{l}$ and $\mathfrak q $ yields $\mathfrak{l}\cdot \mathfrak q =\mathfrak f-f_0+\bar\mathfrak{l}\cdot q_0$. Since $f_0, \bar\mathfrak{l}\cdot q_0 \in \R[x_2,\dots,x_d]$, $\mathfrak{l}(-\bar\mathfrak{l},x_2,\dots,x_d)=\mathfrak f(-\bar\mathfrak{l},x_2,\dots,x_d)=0$ we must have  $-f_0+\bar\mathfrak{l}\cdot q_0=0$.
\end{proof}
Hence the $(d-1)$-dimensional affine subspaces, where a polynomial vanishes are given by the linear factors, and visa versa. Therefore for a set ${\cal I}$ of polynomials in $\R[x_1,\dots,x_d]$ let 
\begin{equation*}
{\cal W}^{d-1}({\cal I})=\{ {\cal C} \subset {\cal V}({\cal I}) :  {\cal C} = \aff({\cal C}) \text{ and } \dim({\cal C})=d-1\}.
\end{equation*}
On account of Remark \ref{linear_factors} we can easily determine  all linear factors and thus ${\cal W}^{d-1}(\mathfrak f)$ of a given polynomial  by the following  procedure.
\begin{algorithm}[Determining a linear factor of a polynomial $\mathfrak f$] \hfill
\label{algorithm_linear_factors}
\begin{itemize} 
\item For each variable $x_j$ do (without loss of generality let $j=1$) 
\begin{enumerate}
\item[{\rm F1.}] Write $\mathfrak f$ as $\mathfrak f=\sum_{i=0}^{m_1} f_{i}\,x_1^i$ with $f_{i}\in \R[x_2,\dots,x_d]$. \par\noindent If $m_1\ne 0$ do 
    \addtolength\itemindent{0.2cm}
\item[{\rm F2.}] Find  $d$ affinely independent points $u^2,\dots,u^{d+1}\in \R^{d-1}$ such that $\mathfrak f(x_1,u^i)\ne 0$, $2\leq i\leq {d+1}$.
\item[{\rm F3.}] For $2\leq i\leq {d+1}$ determine the real roots ${\cal Z}_i$ of the univariate polynomials $\mathfrak f(x_1,u^i)=0$. 
\item[{\rm F4.}] For each $(z_2,\dots,z_{d+1})\in {\cal Z}_2 \times\cdots\times {\cal Z}_{d+1}$  do 
    \addtolength\itemindent{0.2cm}
\item[{\rm F5.}] Determine the  solution $ l=(l_0,l_2,\dots,l_d)^\intercal\in\R^{d}$ of  the linear system  $(1,{u^i}^\intercal) l=z_i$, $2\leq i\leq d+1$.
\item[{\rm F6.}] Let $\mathfrak{l}=-l_0+x_1-\sum_{i=2}^d l_i\cdot x_i$. If $\mathfrak f(x_1-\mathfrak{l},x_2,\dots,x_d)=0$ then $\mathfrak{l}$ is a linear factor and  the remainder $\hat\mathfrak{f}=\mathfrak f/\mathfrak{l}$ is given by the polynomial $\mathfrak q =\sum_{i=0}^{m_1}q_i\, x_1^i$ defined in \eqref{rest}.
\end{enumerate}
 \addtolength\itemindent{0.8cm}
\item If for all solutions $l=(l_0,l_2,\dots,l_d)^\intercal\in\R^{d}$ the corresponding linear polynomial $\mathfrak{l}$ is not a factor of $\mathfrak{f}$ then $\mathfrak{f}$ posseses no linear factors.
\end{itemize}
\end{algorithm} 

\begin{lemma}
\label{algorithm_corect}
Let $\mathfrak f\in\R[x_1,\dots,x_d]$, $\mathfrak f\ne 0$. Then Algorithm {\rm \ref{algorithm_linear_factors}} finds  a linear factor $\mathfrak{l}\in\R[x_1,\dots,x_d]$ of $\mathfrak f$ and  determines the  polynomials $\mathfrak f/\mathfrak{l}$, if $\mathfrak{f}$ has a linear factor at all.  
\end{lemma}
\begin{proof} Let $\mathfrak{l}=l_0+\sum_{i=1}^d l_i x_i$ be a linear factor of  $\mathfrak f$ and let us assume w.l.o.g.~that $l_1=-1$. Let  $\bar\mathfrak{l}=l_0+l_2 x_2+\cdots +l_d x_d\in\R[x_2,\dots,x_d]$ and $l=(l_0,l_2,\dots,l_d)^\intercal$. By Remark \ref{linear_factors} we have $\mathfrak f(\bar\mathfrak{l},x_2,\dots,x_d)=0$ and  for the points $u^i$ we get $0=\mathfrak f(\bar\mathfrak{l}(u^i), u^i)=\mathfrak f((1,{u^i}^\intercal)l, u^i)$ and hence $(1,{u^i}^\intercal)l \in {\cal Z}_i$, $i=2,\dots,d+1$. By the choice of the vectors $u^i$ each linear system has an unique solution.

Finally, we remark that the $(d-1)$-linearly independent points $u^i$  can be found quite easily.  If the leading coefficient of $f$ with respect to $x_1$ is constant, then we set $u^i=e^{i-1}\in\R^{d-1}$, $2\leq i \leq d$, and $u^{d+1}=0$, where $e^i$ denotes the $i$-th unit vector. Otherwise we perturb these points a bit.
\end{proof}

Since each polynomial can be written as an unique product of irreducible polynomials we can apply the above algorithm iteratively to the computed remainders in order to determine all linear factors of a given polynomial. Obviously, we can also use Algorithm \ref{algorithm_linear_factors} to find a common linear factor of two or more polynomials. Hence we  have 

\begin{corollary}
\label{common_factors} Let ${\cal I}$ be as set of polynomials $\mathfrak f_i\in \R[x_1,\dots,x_d]$, $i\in I$. Then there exists an algorithm which computes all common linear factors of the polynomials $\mathfrak f_i$, $i\in I$, and the set ${\cal W}^{d-1}({\cal I})$.
\end{corollary}

For  special polynomials  we have a simple test to decide whether one polynomial is a factor of another one.
\begin{remark} 
\label{test_factor}
Let $\mathfrak{f}=\sum_{i=0}^{m} f_i\cdot x_1^{m-i} \in \R[x_1,\dots,x_d]$, $\mathfrak{g}=\sum_{i=0}^{n} g_i\cdot x_1^{n-i} \in \R[x_1,\dots,x_d]$, $m\geq n\geq 1$, $f_i,g_i\in \R[x_2,\dots,x_d]$  and $g_0\in \R\backslash\{0\}$. Then there exists an algorithm which decides whether $\mathfrak g$ is a factor of $\mathfrak f$ and determines the polynomial $\mathfrak{h}=\mathfrak f /\mathfrak g$  if $\mathfrak g$ is a factor of $\mathfrak f$. 
\end{remark}
\begin{proof} W.l.o.g.~let $g_0=1$. By definition $\mathfrak g$ is a factor of $\mathfrak f$ if and only if there exists a polynomial $\mathfrak h = \sum_{i=0}^{m-n} h_i\cdot x_1^{m-n-i}$, $h_i\in\R[x_2,\dots,x_d]$,  such that $\mathfrak g\cdot  \mathfrak h=\mathfrak f$. Hence by comparing the coefficients we get 
$$ 
   f_k = \sum_{j=\max\{0,n-m+k\}}^{\min\{k,n\}} g_j\cdot h_{k-j}, \quad k=0,\dots.m.
$$
Since $g_0=1$ the first $m-n$ identities determine uniquely  the coefficients $h_j$, $j=0,\dots, m-n$, and the remaining identities can be used as a verification whether $\mathfrak g$ is really a factor of $\mathfrak f$, i.e, whether $\mathfrak h =\mathfrak f /\mathfrak g$.
\end{proof}

Next we describe algorithms how we can find all the isolated affine subspaces for some  very special  polynomials in two variables. 
\begin{lemma}
\label{quadratic_polynomial}
 Let $\mathfrak{f}=f_0\cdot x^2+ f_1\cdot x + f_2\in \R[x,y]$ with $f_i\in \R[y]$ and $\deg(f_i)\leq i$.   Then there exists an algorithm which computes  ${\cal V}_\aff(\mathfrak f)$. 
\end{lemma} 
\begin{proof} Since ${\cal V}(\mathfrak f)$ is a conic section, the set of all isolated affine subspaces of ${\cal V}(\mathfrak f)$ consists   either of  one or of two lines or of an isolated  root or it is empty. Using well know formula for conic sections one can determine all of them.

\end{proof}

\begin{lemma} 
\label{special_polynomial} 
Let $\mathfrak{f}\in\R[x,y]$  with $\deg(\mathfrak f)\leq 4$ and  $\mathfrak{f}\ne 0$. Then there exists an algorithm which computes a set $\widetilde{\cal V}_\aff(\mathfrak f)$ in the sense of \eqref{approxaff}.
\end{lemma}
\begin{proof} Let $\mathfrak{f}=\sum_{i=0}^4 f_i\cdot x^{4-i}$, $f_i\in \R[y]$ and $\deg(f_i)\leq i$.  On account of Corollary \ref{common_factors} we determine all linear factors of $\mathfrak f$ and the set ${\cal W}^1(\mathfrak f)$ containing ${\cal V}^{1}_\aff(\mathfrak f)$.  Hence it remains to determine the isolated roots, i.e., ${\cal V}^0_\aff(\mathfrak f)$. Since none of the isolated roots lies in a $1$-dimensional subspace  of ${\cal W}^1(\mathfrak f)$, we can divide $\mathfrak{f}$ by the linear factors corresponding to the elements of ${\cal W}^1(\mathfrak f)$. 
Therefore, in the following we assume that $\mathfrak{f}$  has no linear factors. If $\mathfrak f$ is an univariate polynomial in $x$, say, having a real root $\bar{x}$, then $-\bar{x}+x$ would be a linear factor of $\mathfrak f$. Thus if $\mathfrak f$ is an univariate polynomial we set $\widetilde{\cal V}_\aff(\mathfrak f)={\cal W}^1(\mathfrak f)$. Hence let $\deg(\mathfrak{f},x) \geq\deg(\mathfrak{f},y)   \geq 1$. 

If $(\bar{x},\bar{y})$ is an isolated root then $\nabla\mathfrak{f}(\bar{x},\bar{y})=0$ and thus we can look for all isolated roots of the set  $\{ (x,y)\in\R^2 :  \mathfrak{f}(x,y)=\mathfrak{f}_x(x,y)=0\}$, 
where $\mathfrak{f}_x$ denotes the partial derivative of $\mathfrak{f}$ with respect to $x$, i.e., 
$
 \mathfrak{f}_x = \sum_{i=0}^3 \tilde{f}_i\cdot x^{3-i},
$
with  $\tilde{f}_i\in \R[y]$, $\deg(\tilde{f}_i)\leq i$. Now let $\res=\res(\mathfrak{f},\mathfrak{f}_x,x)\in \R[y]$ be the resultant of these two polynomials. If $\res \ne 0$ then we determine  ${\cal V}(\res)$ and afterwards ${\cal W}_0=\{ (x,y) : \mathfrak{f}=\mathfrak{f}_x=0,\, y\in {\cal V}(\res)\}$. Observe, that or any fixed $\bar{y}\in {\cal V}(\res)$ the polynomial $\mathfrak{f}(x,\bar{y})$ can not vanish, because otherwise it contains a linear factor (see Remark \ref{linear_factors}). In this case we set $\widetilde{\cal V}_\aff(\mathfrak f)={\cal W}_0\cup {\cal W}^1(\mathfrak f)$.

It remains to consider the case $\res=0$ and therefore we may assume $\deg(\mathfrak{f},x)\geq 2$.  By Lemma \ref{resultants_basic} we know that $\mathfrak{f}$ and $\mathfrak{f}_x$ have a common factor $\mathfrak{g}$ which has positive degree in $x$. Even more, it is not hard to see that $\res(\mathfrak{f},\mathfrak{f}_x)=0$ implies that $\mathfrak{f}$ is divisible by some $\mathfrak{h}^2$, $\mathfrak{h}\in \R[x,y]$ with $\deg(\mathfrak h,x)\geq 1$. Hence, if $\deg(\mathfrak{f},x)\leq 3$ then $\mathfrak{f}$ can be written as a product of polynomials $\mathfrak{f}_k$ with $\deg(\mathfrak{f}_k,x)=1$ which shows that $\mathfrak{f}$ has no isolated roots. 

So let $\deg(\mathfrak{f},x)=4$.  
 Then the common factor $\mathfrak{g}$ of $\mathfrak{f}, \mathfrak{f}_x$ has degree 1, 2 or 3 w.r.t.~the variable $x$. If $\deg(\mathfrak g,x)=1$ than it is a linear factor, because the leading coefficient $f_0$ is a constant. Also, if $\deg(\mathfrak g,x)=3$ than $\mathfrak{f}/\mathfrak{g}\in\R[x,y]$ is a linear factor of $\mathfrak{f}$. Hence $\deg(\mathfrak{g},x)=2$ and it can be written as 
$  \mathfrak{g}=\sum_{i=0}^2 g_i x^{2-i}$   
with $g_i\in \R[y]$, $\deg(g_i)\leq i$. Thus $\mathfrak{f}_x/\mathfrak{g}$ is a linear factor and we can determine $\mathfrak{g}$ by the following procedure: For each linear factor $\mathfrak{l}$ of $\mathfrak{f}_x$ let $\mathfrak{g}_\mathfrak{l}=\mathfrak{f}_x/\mathfrak{l}$. Then we have to test whether $\mathfrak{g}_\mathfrak{l}$ is a factor of $\mathfrak{f}$. This can be done with the algorithm described in Remark \ref{test_factor}.
So we can  assume that we have  found the   factor $\mathfrak{g}$ of $\mathfrak{f}$ with $\deg(\mathfrak{g})=2$.  Let   $\bar{\mathfrak{g}}=\mathfrak{f}/\mathfrak{g}$. Then the problem is reduced to the determination of the isolated roots of the two polynomials $\mathfrak{g}$ and $\bar{\mathfrak{g}}$, which can be solved by Lemma \ref{quadratic_polynomial}. Hence we set $\widetilde{\cal V}_\aff(\mathfrak f)={\cal V}^0_\aff(\mathfrak g)\cup {\cal V}^0_\aff(\bar{\mathfrak g})\cup {\cal W}^1(\mathfrak f)$.
\end{proof}

\begin{lemma} 
\label{two_special_polynomial} 
Let $\mathfrak{f},\mathfrak{g}\in\R[x,y]$ with $\mathfrak{f},\mathfrak{g}\ne 0$ and $\deg(\mathfrak f)\leq 4$, $\deg(\mathfrak{g})\leq 3$. Then there exists an algorithm which computes a set $\widetilde{\cal V}_\aff(\mathfrak f,\mathfrak{g})$ in the sense of \eqref{approxaff}.
\end{lemma}
\begin{proof} We proceed as in the proof of the last lemma. First we determine all common linear factors of these two polynomials and the set ${\cal W}^1(\mathfrak f, \mathfrak g)$ containing ${\cal V}^1_\aff(\mathfrak f, \mathfrak g)$ (cf.~Corollary \ref{common_factors}). Then  we divide the polynomials by these common linear factors and it remains to find the isolated roots ${\cal V}^0_\aff(\mathfrak f, \mathfrak g)$.  To this end we may assume that $\deg(\mathfrak f,x), \deg(\mathfrak f,y) \geq 1$.


Let $\res=\res(\mathfrak f, \mathfrak g,x) \in \R[y]$. If $\res\ne 0$ then we determine ${\cal V}(\res)$, ${\cal W}_0=\{ (x,y) : \mathfrak{f}=\mathfrak{f}_x=0,\, y\in {\cal V}(\res)\}$ and  we set $\widetilde{\cal V}_\aff(\mathfrak f)={\cal W}_0\cup {\cal W}^1(\mathfrak f)$. Therefore we can assume $\res=0$. 
 Since $\deg(\mathfrak g)\leq 3$ we know that if the polynomial $\mathfrak g$ is not irreducible then  it has a linear factor $\mathfrak l$, say. Using Algorithm \ref{algorithm_linear_factors} we can determine such a linear factor $\mathfrak l$ as well as the remainder $\hat\mathfrak g = \mathfrak g/\mathfrak  l$.   Next we split our system in two systems, namely 
\begin{equation*}
 {\cal I}_1  : \mathfrak f = \hat \mathfrak g =0 \text{ and } 
 {\cal I}_2   : \mathfrak f = \mathfrak l =0.
\end{equation*}
Since $\mathfrak l$ is a linear factor and $\mathfrak g$ and $\mathfrak f$ are assumed to be free of common linear factors the second system can be solved by substituting one variable in $\mathfrak f$ via $\mathfrak l$. Since $\deg(\hat\mathfrak g)< \deg(\mathfrak g)$ we can apply 
recursively the previous argumentation to the system ${\cal I}_1$. Thus, if $\mathfrak g$ is not irreducible we set $\widetilde{\cal V}_\aff(\mathfrak f, \mathfrak g)={\cal V}_\aff^0({\cal I}_1) \cup   {\cal V}_\aff^0({\cal I}_2) \cup  {\cal W}^1(\mathfrak f,\mathfrak g)$. 

If $\mathfrak g$ is irreducible, and since $\res =0$, the polynomial $\mathfrak{g}$ has to be a factor of $\mathfrak{f}$. Hence 
 the determination of ${\cal V}^0_\aff(\mathfrak f, \mathfrak g)$ is reduced to the calculation of the isolated roots of $\mathfrak g$. With  the algorithm of Lemma \ref{special_polynomial} we can find a set $\widetilde{{\cal V}}_{\aff}(\mathfrak{g})\subset{\cal V}(\mathfrak{g})$ containing ${\cal V}_{\aff}(\mathfrak{g})$ and we set 
$\widetilde{\cal V}_\aff(\mathfrak f, \mathfrak g)=  \widetilde{{\cal V}}_{\aff}(\mathfrak{g})\cup {\cal W}^1(\mathfrak f, \mathfrak g)$.
\end{proof}

\begin{lemma} 
\label{resultant_affine_spaces}
For $j=1,\dots,3$ let  $\mathfrak{f}^j=f^j_0\cdot x^2 + f^j_1\cdot x + f^j_2 \in \R[x,y,z]$, $f^j_i \in \R[y,z]$, with $\deg(f^j_i)\leq i$.  Let $f^2_0=0$, $f^2_1\ne 0$, $\res_{1,2}=\res(\mathfrak f^1, \mathfrak f^2,x)$, $\res_{2,3}=\res(\mathfrak f^2, \mathfrak f^3,x)$ and let $\res_{1,2}\ne 0$ or $\res_{2,3}\ne 0$. Then for each $L\in {\cal V}_{\aff}(\mathfrak f^1, \mathfrak f^2, \mathfrak f^3)$ there exists a $M\in {\cal V}_{\aff}(\res_{1,2}, \res_{2,3})\cup W^1(\res_{1,2},\res_{2,3})$ such that $L\subset \{(x,y,z) \in \R^3 : (y,z)\in M\}$. 
\end{lemma}
\begin{proof} W.l.o.g.~let $\res_{1,2}\ne 0$ and let $L\in {\cal V}_{\aff}(\mathfrak f^1, \mathfrak f^2, \mathfrak f^3)$. By $L^\pi$ we denote the orthogonal projection of $L$ onto the plane $\{ (x,y,z)\in \R^3 : x=0\}$. Obviously, $L^\pi$ is an affine subspace and by the definition of resultants we have $\res_{1,2}(y,z)=\res_{2,3}(y,z)=0$ for all $(y,z)\in L^\pi$. Since $\res_{1,2}\ne 0$ we have $\dim(L^\pi)\in \{0,1\}$. If $\dim(L^\pi)=1$, then $L^\pi$ itself corresponds to a common linear factor of $\res_{1,2}$ and $\res_{2,3}$ and thus $L^\pi\in W^1(\res_{1,2},\res_{2,3})$. Therefore let $L^\pi=\{ (y^0,z^0)\}$ and so $L$ is either an isolated $1$-dimensional subspace or an isolated root of ${\cal V}(\mathfrak f^1, \mathfrak f^2, \mathfrak f^3)$. If $(y^0,z^0)$ is contained in an $1$-dimensional subspace of $W^1(\res_{1,2}, \res_{2,3})$ the   statement is certainly true. Thus we may assume that $\mathfrak l(y^0,z^0)\ne 0$ for every common linear factor $\mathfrak l$ of $\res_{1,2}$ and $\res_{2,3}$ and we have to show that $u^0=(y^0,z^0)$ is an isolated root of ${\cal V}_\aff(\res_{1,2}, \res_{2,3})$. Suppose the contrary and let $u^1=(y^1,z^1) \in {\cal V}(\res_{1,2}, \res_{2,3})$ such that there exists a path ${\cal P} \subset {\cal V}(\res_{1,2}, \res_{2,3})$, ${\cal P}=\{ u^t : t\in[0,1]\}$, connecting $u^0$ and $u^1$. Then ${\cal P} \ne \conv\{u^0,u^1\}$, because otherwise $\aff\{u^0,u^1\}\subset {\cal V}(\res_{1,2}, \res_{2,3})$ is  a $1$-dimensional set containing $u^0$. Hence we can assume $f^2_1(u^t)\ne 0$, $t\in(0,1)$, and by Lemma \ref{resultants_basic} ii) there exist $a^t,b^t \in \C$ such that $\mathfrak f^1(a^t,u^t)=\mathfrak f^2(a^t,u^t)= 0 = f^2(b^t,u^t) = \mathfrak f^3(b^t,u^t)$, $t\in (0,1)$. Since $f^2_0=0$ and  $f^2_1(u^t)\ne 0$ we get $a^t=b^t\in \R$. However this shows that $L$ is not an isolated affine subspace of ${\cal V}(\mathfrak f^1, \mathfrak f^2, \mathfrak f^3)$. 
\end{proof}
Of course the last lemma also implies 
\begin{corollary} 
\label{resultant_affine_spaces_cor}
For $j=1,2$ let  $\mathfrak{f}^j=f^j_0\cdot x^2 + f^j_1\cdot x + f^j_2 \in \R[x,y,z]$, $f^j_i \in \R[y,z]$, with $\deg(f^j_i)\leq i$.  Let $f^2_0=0$, $f^2_1\ne 0$ and let  $\res_{1,2}=\res(\mathfrak f^1, \mathfrak f^2,x)\ne 0$. Then for each $L\in {\cal V}_{\aff}(\mathfrak f^1, \mathfrak f^2)$ there exists a $M\in {\cal V}_{\aff}(\res_{1,2})\cup W^1(\res_{1,2})$ such that $L\subset \{(x,y,z) \in \R^3 : (y,z)\in M\}$. 
\end{corollary}
\begin{proof} We set $f^3=f^1$ and apply Lemma \ref{resultant_affine_spaces}.
\end{proof}

Using Corollary \ref{common_factors} and the previous lemmas  we can make the resultant approach practicable. 
\begin{theorem} 
\label{polynomial_theorem}
Let $\mathfrak{p}(x,y,z)$ be a polynomial with total degree at most 3. Let ${\cal V}(\nabla \mathfrak p)=\{ (x,y,z)\in \R^3 : \nabla p =0\}$. There exists an algorithm 
 computing a set $\widetilde{\cal V}_\aff(\nabla \mathfrak p)$ in the sense of \eqref{approxaff}. \end{theorem}
\begin{proof} First we note that after some scaling, subtractions, and renumbering we may assume that ${\cal V}(\nabla \mathfrak{p})$ is given by  the following system ${\cal I}$ of polynomial equations (see \eqref{partial})
\begin{equation}
\label{partial2}
{\cal I} : \quad\quad 
\begin{split}
 \mathfrak{p}_1 & = \hphantom{ \,\kappa\cdot x^2 +}\,    \mathfrak{l}_1\cdot x  + 
                  \mathfrak{q}_1 =0 , \\
 \mathfrak{p}_2 & = \hphantom{ \,\kappa\cdot x^2 +}\,   \mathfrak{l}_2\cdot x + 
                  \mathfrak{q}_2 =0, \\
 \mathfrak{p}_3 & =  \kappa\cdot x^2 + \mathfrak{l}_3\cdot x + 
                  \mathfrak{q}_3 =0,
\end{split}
\end{equation}
where $\kappa\in\{0,1\}$, $\mathfrak{l}_i\in\R[y,z]$ are linear polynomials and $\mathfrak{q}_i\in\R[y,z]$ with $\deg(\mathfrak q_i)\leq 2$. We may further assume that $\deg(\mathfrak{p}_1)\leq \deg(\mathfrak{p}_2) \leq \deg(\mathfrak{p}_3) $.  Depending on the number of non-trivial polynomials and the number of variables in ${\cal I}$ we have to distinguish several cases. Obviously, if all polynomials in \ref{partial2} vanish then we have ${\cal V}(\nabla \mathfrak p)=\R^3$ and we have to do nothing. 
\begin{enumerate}

\item[(0)] ${\cal I}$ consists of one or two or three polynomials in only one variable. 
\par\noindent Then we can determine ${\cal V}(\nabla \mathfrak p)$ by any algorithm computing  the roots of an univariate polynomial. 

\item[(1,2)] ${\cal I}$ consists of one polynomial in two variables.
\par\noindent W.l.o.g.~let  ${\cal V}({\cal I})=\{(x,y,z)\in \R^3 : \mathfrak{q}_3(y,z)=0\}= \R \times {\cal V}(\mathfrak q_3)$, with ${\cal V}(q_3)\subset \R^2$.  Via the algorithm of Lemma \ref{quadratic_polynomial} we can determine ${\cal V}(\mathfrak q_3)$ and we set $\widetilde{\cal V}_\aff({\cal I})=\{ \R\times {\cal C} : {\cal C} \in   {\cal V}_\aff(\mathfrak q_3)\}$.

\item[(2,2)] ${\cal I}$ consists of two polynomials in two variables.
\par\noindent W.l.o.g.~let ${\cal V}({\cal I})=\{(x,y,z)\in \R^3 : \mathfrak{q}_2(y,z)=\mathfrak{q}_3(y,z)=0\}=\R \times {\cal V}(\mathfrak q_2, \mathfrak q_3)$. Using the Algorithm of Lemma \ref{two_special_polynomial} we can determine a set $\widetilde{\cal V}_\aff(\mathfrak q_2,\mathfrak q_3)\subset \R^2$ (cf.~\eqref{approxaff}) and  we set $\widetilde{\cal V}_\aff({\cal I})=\{ \R\times {\cal C} : {\cal C} \in   \widetilde{\cal V}_\aff(\mathfrak q_2,\mathfrak q_3)\}$.

\item[(3,2)] ${\cal I}$ consists of three polynomials in two variables.
\par\noindent W.l.o.g.~let ${\cal V}({\cal I})=\{(x,y,z)\in \R^3 : \mathfrak{q}_1(y,z)=\mathfrak{q}_2(y,z)=\mathfrak{q}_3(y,z)=0\}=\R \times {\cal V}(\mathfrak q_1,\mathfrak q_2, \mathfrak q_3)$. We may assume that both variables occur in all three  polynomials and that the polynomials are linearly independent. Otherwise we can reduce  this case to one of the previous ones.  First, by Corollary \ref{common_factors}  we determine all common linear factors and the set ${\cal W}^1(\mathfrak q_1,\mathfrak q_2, \mathfrak q_3)$ containing ${\cal V}^1_\aff(\mathfrak q_1,\mathfrak q_2, \mathfrak q_3)$ and hence we may assume that   $\mathfrak q_1, \mathfrak q_2, \mathfrak q_3$ have no common linear factors and  both variables occur in the polynomials. Next we compute the resultant $\res_{2,3}=\res(\mathfrak{q}_2,\mathfrak{q}_3,y)\in \R[z]$. If $\res_{2,3}\ne 0$ then we determine ${\cal W}=\{ (y,z) : \mathfrak q_1 = \mathfrak q_2 =\mathfrak q_3 = 0, \, z \in {\cal V}(\res_{2,3})\}$. Observe that for a fixed $\bar{z}$ not all three polynomials $\mathfrak q_i$ can vanish, because  otherwise they have a common linear factor. In this case we set $\widetilde{\cal V}_\aff({\cal I})= \{     \R\times {\cal C} : {\cal C} \in {\cal W}^1(\mathfrak q_1,\mathfrak q_2, \mathfrak q_3) \cup {\cal W}\}$. \newline    
It remains to consider the case $\res_{2,3}=0$. Since $\mathfrak{q}_2$, $\mathfrak{q}_3$ are assumed to be linearly independent and since $\deg(\mathfrak q_i)\leq 2$  the common factor (cf.~Lemma \ref{resultants_basic} iii)) has to be  a linear polynomial $\mathfrak l_{2,3}$ which can be determined  via  Algorithm \ref{algorithm_linear_factors}. Then we consider the two systems 
\begin{align*}
{\cal I}_1 & :\quad \mathfrak{q_1}=0,\,  \mathfrak{l}_{2,3}=0 \\ 
{\cal I}_2 & :\quad \mathfrak{q_1}=0, \, \mathfrak{q}_2/\mathfrak{l}_{2,3}=0, \,\mathfrak{q}_3/\mathfrak{l}_{2,3}=0.
\end{align*}
Both systems are free of common linear factors and since both systems contain linear polynomials  we can easily determine ${\cal V}({\cal I}_1)$ and  ${\cal V}({\cal I}_2)$ by substitution. We set  $\widetilde{\cal V}_\aff({\cal I})= \{ \R\times {\cal C} : {\cal C} \in {\cal W}^1(\mathfrak q_1,\mathfrak q_2, \mathfrak q_3) \cup {\cal V}({\cal I}_1) \cup {\cal V}({\cal I}_2) \}$.

\item[$\bullet$] In the remaining cases we first determine the set ${\cal W}^2({\cal I})$ and therefore we may always assume that the given partial derivatives have no common linear factors.

\item[(1,3)] ${\cal I}$ consists of one  polynomial in three  variables.
\par\noindent {\rm a )} Let ${\cal V}({\cal I})=\{(x,y,z) \in \R^3 : \mathfrak{l}_3\cdot x + \mathfrak{q}_3 = 0\}$ Then we may assume $\mathfrak{l}_3\ne 0$ and the dimension of any  affine subspace of ${\cal V}({\cal I})$ is 2. Therefore we set $\widetilde{\cal V}({\cal I})={\cal W}^2({\cal I})$.
\par\noindent {\rm b )} Let ${\cal V}({\cal I})=\{(x,y,z) \in \R^3 : \mathfrak f = x^2+ \mathfrak{l}_3\cdot x + \mathfrak{q}_3 = 0\}$. 
Let $\tilde\mathfrak{q}=\mathfrak{l}_3\cdot \mathfrak{l}_3/4-\mathfrak{q}_3\in \R[y,z]$. For every $(x^*,y^*,z^*)\in{\cal V}({\cal I})$ we have $x^*=-\mathfrak{l}_3(y^*,z^*)/2\pm \sqrt{ \tilde\mathfrak{q}(y^*,z^*)}$ and hence  $\tilde\mathfrak{q}(y^*,z^*)\geq 0$. If $\tilde\mathfrak{q}(y^*,z^*)>0$ 
 then we can always find  a  neighbourhood $U$ of $(y^*,z^*)$ such that $\tilde\mathfrak{q}(y,z)\geq 0$ for all $(y,z)\in U$ and  $(x^*,y^*,z^*)$ belongs  to a 2-dimensional connected component of ${\cal V}(\mathfrak{f})$. Hence we may set (cf.~Lemma \ref{quadratic_polynomial})
$$
\widetilde{\cal V}_\aff({\cal I})={\cal W}^2({\cal I})\cup_{{\cal C}\in {\cal V}_\aff(\tilde\mathfrak{q})}
 \left\{ (-\mathfrak{l}_3(y,z)/2,y,z ) : (y,z)\in {\cal C} \right\}.
$$

\item[(2,3)] ${\cal I}$ consists of two  polynomials in three  variables.
\par\noindent Then we may assume w.l.o.g.~that the variable $x$ occurs in both polynomials.  
\par\noindent {\rm a )} Let ${\cal V}({\cal I})=\{(x,y,z) \in \R^3 : \mathfrak p_1=\mathfrak{l}_1\cdot x + \mathfrak{q}_1 = 0, \,\mathfrak p_2=\mathfrak{l}_2\cdot x + \mathfrak{q}_2 = 0\}$ with  $\mathfrak{l}_1\ne 0$, $\mathfrak{l}_2\ne 0$ and  we can assume that $\mathfrak p_1, \mathfrak p_2$ are linearly independent. 
Let $\res_{1,2}\in\R[y,z]$ be the resultant of these two polynomials with respect to $x$. Since the polynomials are linearly independent, $\deg(\mathfrak p_1),\deg(\mathfrak p_2)\leq 2$, and since we have assumed that they have no common linear factors  the resultant $\res_{1,2}$ can not vanish.   By definition $\res_{1,2}\in \R[y,z]$ is a polynomial of total degree at most 3. By Lemma \ref{special_polynomial} we can find a set $\widetilde{\cal V}_\aff(\res_{1,2})$ containing ${\cal V}_\aff(\res_{1,2})$ and for each ${\cal C}\in \widetilde{\cal V}_\aff(\res_{1,2})$ we consider the system 
$$
  {\cal I}_{\cal C} : \mathfrak p_1(x,y,z) =0,\, \mathfrak p_2(x,y,z)=0, \, (y,z)\in{\cal C}.
$$
Since ${\cal C}$ is a $0$- or $1$-dimensional affine subspace, ${\cal I}_{\cal C}$ is a system of at most two polynomials in at most two variables. Hence by the previous cases we can determine a set $\widetilde {\cal V}_\aff({\cal I}_{\cal C})\subset  
{\cal V}({\cal I}_{\cal C})$ containing $ {\cal V}_\aff({\cal I}_{\cal C})$ and we set 
$$
\widetilde{\cal V}_\aff({\cal I})={\cal W}^2({\cal I})\bigcup_{{\cal C}\in \widetilde{\cal V}_\aff(\res_{1,2})}  \widetilde{\cal V}_\aff({\cal I}_{\cal C}).
$$
On account of Corollary \ref{resultant_affine_spaces_cor} we have ${\cal V}_\aff({\cal I})\subset \widetilde{\cal V}_\aff({\cal I})$.

\par\noindent {\rm b )} Let ${\cal V}({\cal I})=\{(x,y,z) \in \R^3 : \mathfrak p_1=\mathfrak{l}_1\cdot x + \mathfrak{q}_1 = 0, \,\mathfrak p_2=x^2+\mathfrak{l}_2\cdot x + \mathfrak{q}_2 = 0\}$ with $\mathfrak l_1 \ne 0$. Then we can proceed as in the case above. The only difference is that the total degree of $\res_{1,2}$ is at most 4.

\item[(3,3)] ${\cal I}$ consists of three  polynomials in three  variables.
\par\noindent Then we may assume w.l.o.g.~that the variable $x$ occurs in at least two of the  polynomials,  $\mathfrak l_2\ne 0$ (cf.~\eqref{partial2}) and that each of the three polynomials contains at least two variables. 
Now we compute $\res_{1,2}=\res(\mathfrak p_1,\mathfrak p_2,x)$ and $\res_{2,3}=\res(\mathfrak p_2,\mathfrak p_3,x)$. 
\par\noindent a) $\res_{1,2}=0$. 
\par\noindent By Lemma \ref{resultants_basic} we know that  $\mathfrak p_1$ and $\mathfrak p_2$ have a common factor with positive degree in $x$. Thus $\mathfrak p_1$, $\mathfrak p_2$ are either linearly dependent or they have a common linear factor. If they are linearly dependent we can proceed as in case (2,3) b). So let $\mathfrak l$ be a common linear factor and let $\tilde\mathfrak p_1=\mathfrak p_1/\mathfrak l$, 
$\tilde\mathfrak p_2=\mathfrak p_2/\mathfrak l$. Observe,  $\tilde\mathfrak p_1,\tilde\mathfrak p_2$ are linear polynomials. Next we consider the two systems 
\begin{equation*}
\begin{split}
 {\cal I}_1 & : \mathfrak l=0,\,\, \mathfrak p_3=0, \\
 {\cal I}_2 & : \tilde\mathfrak p_1=0,\,\,\tilde\mathfrak p_2=0,\,\, \mathfrak p_3=0.
\end{split}
\end{equation*}
Since ${\cal V}_\aff({\cal I})\subset {\cal V}_\aff({\cal I}_1)\cup {\cal V}_\aff({\cal I}_2)$ it suffices to consider ${\cal V}_\aff({\cal I}_i)$. Since both systems contain linear polynomials we can reduce them to systems in at most 2 variables which can be handled by one of the methods described in one of the previous cases.

\par\noindent b) $\res_{2,3}=0$. 
It is not hard to see that also in this case  we can split our system in two systems containing linear polynomials and we can proceed as before.

\par\noindent c) $\res_{1,2}\ne 0$ and $\res_{2,3}\ne 0$. 
Since $\deg(\res_{1,2})\leq 3$ and $\deg(\res_{2,3})\leq 4$ we can use the algorithm of Lemma \ref{two_special_polynomial} in order to determine a set $\widetilde{\cal V}_\aff(\res_{1,2},\res_{2,3})$ containing ${\cal V}_\aff(\res_{1,2},\res_{2,3})$. Now, for each ${\cal C}
\in  \widetilde{\cal V}_\aff(\res_{1,2},\res_{2,3})$ let 
$$
 {\cal I}_{\cal C} : \mathfrak p_i(x,y,z)=0,\, i=1,\dots,3,\, (y,z)\in{\cal C}. 
$$
Since ${\cal C}$ is a $0$- or $1$-dimensional affine subspace ${\cal I}_{\cal C}$ is a system of at most three polynomials in at most two variables. Hence by the previous cases we can determine a set $\widetilde {\cal V}_\aff({\cal I}_{\cal C})$ containing 
${\cal V}_\aff({\cal I}_{\cal C})$ and we set 
$$
\widetilde{\cal V}_\aff({\cal I})={\cal W}^2({\cal I})\bigcup_{{\cal C}\in \widetilde{\cal V}_\aff(\res_{1,2},\res_{2,3})}  \widetilde{\cal V}({\cal I}_{\cal C}).
$$
On account of Lemma \ref{resultant_affine_spaces} we have ${\cal V}_\aff({\cal I})\subset \widetilde{\cal V}_\aff({\cal I})$. 
\end{enumerate}
\end{proof}

\par\noindent \textbf{Remark.}
Up to now we have only discussed polynomials arising in case I.~of  Algorithm \ref{algorithm_2}. Case II.~can be treated completely similar to case I. In the cases III.~and IV.~we have seven hyperplanes determining the function $f_{H_{l_1},\dots,H_{l_7}}$, but it is easy to see that Lemma \ref{rang_lemma} keeps true. Hence also in these cases we just  have to examine the isolated affine subspaces of polynomials of total degree at most three. Indeed, with some extra effort one can show that in the cases III.~and IV.~it suffices to consider matrices $(A_{i_1,\dots,i_7})$ with  ${\rm rank}(A_{i_1,\dots,i_7})=7$. However, we note that in case IV.~we have to replace the $6$-th hyperplane $H_{l_6}$ by the hyperplane $\tilde H_{l_6}$  given in Theorem \ref{ness_con}. The hyperplane $\tilde H_{l_6}$ can easily be constructed, if it exists at all. 

Altogether, using the notation of Algorithm \ref{algorithm_2}  we have the following result.

\begin{corollary}
\label{main_corollary}
 There exists an algorithm which computes a  set ${\cal V}_{H_{l_1},\dots,H_{l_k}}$ consisting of finitely many affine subspaces $A_1,\dots, A_r$ of ${\cal S}_{H_{l_1},\dots,H_{l_k}}$ such that $f_{H_{l_1},\dots,H_{l_k}}$ is constant on any affine subspace $A_i$ and each local minimum of $f_{H_{l_1},\dots,H_{l_k}}$ is contained in one of the spaces $A_i$. 
\end{corollary}

Observe,  we do not determine  the local minima ${\cal M}_{H_{l_1},\dots,H_{l_k}}$ of the function $f_{H_{l_1},\dots,H_{l_k}}$, but just a set of affine subspaces containing this set.  
It remains to check whether such an affine subspace $A_i$ of the corollary contains an admissible lattice. To this end we use the criteria of Minkowski as given in Lemma \ref{lattice_check}, i.e., we are only interested in admissible lattices having a basis $W\in A_i$ such that  the conditions of Lemma \ref{lattice_check} are satisfied. 
\begin{lemma} 
\label{lemma_solution}
There exists an algorithm which finds for any affine subspace $A \in {\cal V}_{H_{l_1},\dots,H_{l_k}}$ a basis $W$ of an admissible lattice (in the sense of Lemma {\rm \ref{lattice_check}}) or asserts that no such lattice exists.   
\end{lemma}
\begin{proof}
For simplification we assume that $H_{l_i}=H_i$, $1\leq i\leq k$. 
With  $r=\dim A$, $r\in\{0,3\}$, we may write 
$$
  A=\left\{ W\in\R^{3\times 3} : W = C + \sum_{i=1}^r \lambda_i M_i\right\},
$$
for suitable matrices $C,M_i\in \R^{3\times 3}$. Now for $r\geq 1$ and for  $\lambda=(\lambda_1,\dots,\lambda_r)\in\R^r$ let $W(\lambda)=C + \sum_{i=1}^r \lambda_i M_i$. If $r=0$ we set $W(\lambda)=C$. According to Lemma \ref{lattice_check} and Theorem \ref{ness_con} we have to find  a $\lambda \in\R^r$ such that 
\begin{equation}
\label{solution}
    {\cal U}^j_{W(\lambda)} \subset {\rm bd} P,
\end{equation}
with $j=1$ in case I., $j=2$ in case II., and  $j=3$ in the cases III.~and IV. Furthermore in the cases I.~and II.~we have the additional restrictions 
\begin{equation}
\label{additional_restrictions}
\begin{split}
& \text{Case \,\,I: }   (-1,1,1)_{W(\lambda)},\, (1,-1,1)_{W(\lambda)},\,(1,1,-1)_{W(\lambda)} \notin P, \\
& \text{Case II: }  (1,1,1)_{W(\lambda)}\notin P.
\end{split}
\end{equation}
Let us denote by $u^i_{W(\lambda)}$, $1\leq i\leq k$, the vectors of the test set ${\cal U}^j_{W(\lambda)}$. This means we have $k=6$ if $j\leq 2$, and $k=7$ otherwise. Observe, by construction we know that all vectors $u^i_{W(\lambda)}$ lie in supporting hyperplanes of the polytope. Therefore, if $r=0$  we can verify \eqref{solution} and \eqref{additional_restrictions} by just checking the  facet defining inequalities of the polytope for the corresponding points $u^i_{W(\lambda)}$, etc. Thus in the following we assume $r\geq 1$. Then we define a  set ${\cal A}$ by  (cf.~notation at the beginning of this section)
\begin{equation}
\label{lp_solution}
\begin{split} 
 {\cal A} & =\left\{\lambda\in\R^r : {\cal U}^j_{W(\lambda)} \subset \bd(P)\right\} \\ & = 
 \big\{\lambda \in \R^r : u^i_{W(\lambda)} \in H^+_{i_j}  \text{ for all } 
                  F_{i_j}\in{\cal N}(F^i),\;1\leq i \leq k \big\}.
\end{split}
\end{equation} 
Using standard methods from Linear Programming we can easily decide whether ${\cal A}=\emptyset$ or we can find a point $\lambda^*\in {\cal A}$.  If ${\cal A}=\emptyset$ then the affine subspace does not contain an admissible lattice. So let  $\lambda^*\in {\cal A}$. In the cases III.~and IV.~or if $W(\lambda^*)$ also satisfies \eqref{additional_restrictions} in the cases I.~and II., we are done and we have found an admissible lattice. So let us assume that we are in case I.~or II., ${\cal A}\ne\emptyset$ and $W(\lambda^*)$ violates \eqref{additional_restrictions}. 
The most simple  way to decide whether there exists a $\lambda\in{\cal A}$ satisfying \eqref{additional_restrictions} is the following: In case I. we consider for $i,j,l\in\{1,\dots,n\}$ the sets
\begin{equation} 
\label{caseI_lp}
\begin{split}
  {\cal A}_{i,j,l}=\big\{ \lambda\in\R^r : \, & \lambda\in{\cal A}\quad\text{ and } 
                    (-1,1,1)_{W(\lambda)} \in H^-_i, \,  \\
                  &  (1,-1,1)_{W(\lambda)} \in H^-_j, \, 
                    (1,1,-1)_{W(\lambda)} \in H^-_l 
                   \big\}.
\end{split}
\end{equation}
In case II~we set for $1\leq i\leq n$ 
\begin{equation} 
\label{caseII_lp} 
  {\cal A}_{i}=\big\{ \lambda\in\R^r :  \lambda\in{\cal A}\quad\text{ and } 
                    (1,1,1)_{W(\lambda)} \in H^-_i
                   \big\}.
\end{equation}
Obviously, in the case I.~(case II.) there exists an admissible lattice in the affine subspace $A$ if and only if there exist $i,j,l\in\{1,\dots,n\}$  ($i\in\{1,\dots,n\}$) such that ${\cal A}_{i,j,l}\ne\emptyset$ (${\cal A}_{i}\ne\emptyset$). Again, using tools from Linear Programming we can either find a $\lambda^*$ in one of sets ${\cal A}_{i,j,l}$ (${\cal A}_i$), and thus an admissible lattice $W(\lambda^*)$, or we know that $A$ contains no admissible lattices. 
\end{proof}

\par\noindent \textbf{Remark.} Instead of considering  the $n^3$ ($n$) feasibility problems in \eqref{caseI_lp} (\eqref{caseII_lp}) one can argue as follows: 
First let us assume that we are in case I.~and  that  $W(\lambda^*)$ violates one of the restrictions in \eqref{additional_restrictions}. Then we claim that we do not have to consider this case further: Of course, if each $W(\lambda)$ violates \eqref{additional_restrictions} for all $\lambda\in {\cal A}$ satisfying \eqref{solution} then this is trivial. Hence suppose that there exists a $\lambda^\circ\in{\cal A}$ satisfying \eqref{additional_restrictions} such that $W(\lambda^\circ)\Z^3$ is  a critical lattice of $P$. Then there exists a $\mu \in \conv\{\lambda^\circ,\lambda^*\}\subset {\cal A}$ such that one of the points of \eqref{additional_restrictions} w.r.t.~$W(\mu)$ lies in the boundary of $P$ and the other points are not contained in the interior of $P$. By definition, $W(\mu)\Z^3$ is a critical lattice of $P$, too.
W.l.o.g.~let $(1,1,-1)_{W(\mu)}\in {\rm bd}P$ and  $(1,-1,1)_{W(\mu)}, (-1,1,1)_{W(\mu)}\notin \inte P$. Furthermore, let $H_7$ be a supporting hyperplane of a facet of $P$ containing $(1,1,-1)_{W(\mu)}$. Now let $w^i(\mu)$ be the $i$-th column vector of $W(\mu)$ and let $\bar W(\mu)$ be the matrix with columns $w^1(\mu)-w^3(\mu)$, $w^2(\mu)$ and $-w^2(\mu)+w^3(\mu)$. Then $\bar W(\mu)$ is just another basis of the lattice $W(\mu)\Z^3$ and we get 
\begin{alignat*}{2}
 (1,0,0)_{W(\mu)}&= (1,1,1)_{\bar W(\mu)}, \, 
 (0,1,0)_{W(\mu)}&= (0,1,0)_{\bar W(\mu)}, \, \\
 (0,0,1)_{W(\mu)}&= (0,1,1)_{\bar W(\mu)}, \, 
(1,-1,0)_{W(\mu)}&= (1,0,1)_{\bar W(\mu)}, \, \\
(1,0,-1)_{W(\mu)}&= (1,0,0)_{\bar W(\mu)}, \,
(0,-1,1)_{W(\mu)}&= (0,0,1)_{\bar W(\mu)}, \, \\
(1,1,-1)_{W(\mu)}&= (1,1,0)_{\bar W(\mu)}.
\end{alignat*}
Thus the test set ${\cal U}^1_{W(\mu)}$ plus the additional point $(1,1,-1)_{W(\mu)}$ is equivalent to the test set ${\cal U}^3_{\bar W(\mu)}$. Since $\bar W(\mu)\Z^3$ is a critical lattice of $P$ it follows from the work of Minkowski (cf.~\cite{m}, pp.27) that we shall find a basis of this lattice in case III. 
In case II.~we can apply the argumentation and this means that in all cases it is sufficient to determine only one point of the set ${\cal A}$.

\vspace{0.2cm}
To sum it up, we finally have the following Algorithm for determining a densest lattice packing of a 3-dimensional polytope $P\in {\cal K}^3$ (cf.~Algorithm \ref{algorithm_2}):

\begin{algorithm}[Densest lattice packing of a 3-polytope] \hfill
\label{algorithm_final}
\par\noindent {\rm Input:} A polytope $P\in {\cal K}^3$ given by the supporting hyperplanes $H_i$, $1\leq i\leq m$, or by its vertices. \hfill 
\par\noindent {\rm Output:} A densest packing lattice of $P$.
\vspace{-0.2cm}
\begin{itemize}
\item Find the supporting hyperplanes $H_i$, $1\leq i\leq n$, of the facets $F_i$ of the polytope $P_0=(P-P)\in{\cal K}_0^3$ {\rm (cf.~\eqref{sym_body})} and compute the lattice description of $P_0$. With respect to $P_0$ do
\item For each of the  cases {\rm I., II., III., IV.}~of Theorem {\rm \ref{ness_con}} do \vspace{-0.2cm}
\begin{itemize}
\item [{\rm I1.}] Compute the sets ${\cal G}(F_i)$ and  ${\cal G}$ with Algorithm {\rm \ref{algorithm_g}}.
\item For three facets  $F_{l_1},F_{l_2},F_{l_3}$   satisfying $F_{l_2} \in{\cal G}(F_{l_1})$ and $F_{l_3}\in {\cal G}(F_{l_1})\cap {\cal G}(F_{l_2})$ do  
    \addtolength\itemindent{0.2cm}

\item For every choice of  facets $F_{l_i}$, $4\leq i\leq k$,  with $(F_{l_1},F_{l_2},F_{l_6}), \allowbreak  (F_{l_2},F_{l_3},F_{l_4}),  (F_{l_1},F_{l_3},F_{l_5}) \in {\cal G}$ do 
\begin{enumerate}
    \item[{\rm R0.}] If ${\widetilde {\cal S}}_{H_{l_1},\dots,H_{l_k}}\ne \emptyset$ {\rm (cf.~\eqref{feaslp})} do 
    \addtolength\itemindent{0.2cm}
    \item[{\rm R1.}] Determine ${\cal S}_{H_{l_1},\dots, H_{l_k}}$ {\rm (cf.~\eqref{para_set})}. \newline 
If ${\rm rank}(A_{i_1,\dots,i_k})\geq 6$ {\rm (cf.~Lemma \ref{rang_lemma} and the remark after Theorem \ref{polynomial_theorem}) } do 
    \addtolength\itemindent{0.1cm}
    \item[{\rm R2.}]  Determine a set ${\cal V}_{H_{l_1},\dots,H_{l_k}}$ consisting of finitely many affine subspaces $A_1,\dots, A_r$ of ${\cal S}_{H_{l_1},\dots,H_{l_k}}$ such that $f_{H_{l_1},\dots,H_{l_k}}$ is constant on any affine subspace $A_i$ and each local minimum of $f_{H_{l_1},\dots,H_{l_k}}$ is contained in one of the spaces $A_i$ {\rm (cf.~Corollary \ref{main_corollary})}. 
\item[{\rm R3.}] For each affine  subspace $A \in {\cal V}_{H_{l_1},\dots,H_{l_k}}$ find a basis $W$ such that $W\cdot \Z^d$ satisfies the criterion {\rm 1}.~of Lemma {\rm \ref{lattice_check}} in the first case and criterion {\rm 2.}~in the remaining cases,  or asserts that  such lattice does not  exist {\rm (cf.~Lemma \ref{lemma_solution})}.   
\end{enumerate}
\end{itemize}
\item Among all calculated admissible lattices find one with  minimal determinant. The corresponding lattice is a critical lattice of $P_0$ and a densest packing lattice of $P$ .
\end{itemize}
\end{algorithm}

\section{Densest lattice packings of the regular and \\  Archimedean  polytopes}

In this section we present densest packing lattices $\Lambda^*$ of all regular and Archimedean polytopes. Since a densest packing lattice depends on the representation of the polytope, and in order to make the results more transparent, we shall also give the coordinates of our representations of the polytopes used for the algorithm. For a polytope $P$ let $f=(f_0,f_1,f_2)$ be its $f$-vector, i.e., $f_i$ is the number of $i$-faces. Futhermore we shall use the following abbrevations: 
 Let  $\tau=(1/2)(1+\sqrt{5})$ and let $D_3$ be the fcc-lattice with basis $(1,1,0)^\intercal$, $(1,0,1)^\intercal$, $(0,1,1)^\intercal$. Finally, let $P_{rtc}$ be the so called rhombic triacontahedron given by 
\begin{equation}
\begin{split}
 P_{rtc}=\Big\{ x\in \R^3 : & \,|\tau x_i| \leq 1, \, \left|\frac{1}{2}x_1 + \frac{\tau}{2}x_2 + \frac{\tau+1}{2}x_3 \right| \leq 1,\, \\ 
&\left|\frac{\tau}{2}x_1 + \frac{\tau+1}{2}x_2 + \frac{1}{2}x_3 \right| \leq 1,\,
\left|\frac{\tau+1}{2}x_1 + \frac{1}{2}x_2 + \frac{\tau}{2}x_3 \right| \leq 1\Big\}. 
\end{split}
\end{equation}
For the identification of the polytopes we shall use the Wythoff symbols. 
\ifnum\pictures=1 
 Pictures of the polytopes are given at the end of this section.
\fi 

\begin{itemize}

\item {\sl Tetrahedron}, $3\,|\,2\,3$, $f=(4,6,4)$ (cf.~\cite{h}).
\begin{equation*}
\begin{split}
  &         P_{t}    = \big\{x \in \R^3 : x_1+x_2+x_3\leq 1, \, 
                                -x_1-x_2+x_3\leq 1,\, \\
  & \hphantom{P_{t}    = \big\{x \in \R }  
        -x_1+x_2-x_3\leq 1,\, \,\,\,\, x_1-x_2-x_3\leq 1 
                             \big\}, \\ 
 &  \Lambda^*(P_{t})  = 2\left( (1, -\frac{1}{6}, -\frac{1}{6})^\intercal,
                (-\frac{1}{6}, 1, -\frac{1}{6})^\intercal,
                (-\frac{1}{6}, -\frac{1}{6}, 1)^\intercal 
                 \right)\Z^3, \\
 & \delta^*(P_{t})   = \frac{18}{49}\approx 0.367346938.
\end{split}
\end{equation*}
\item {\sl Cube}, $3\,|\,2\,4$, $f=(8,12,6)$.
\begin{equation*}
\begin{split} 
&  P_{c}    = \left\{x \in \R^3 : |x_i| \leq 1 \right\}, \quad                 
  \Lambda^*(P_{c})  = 2\Z^3, \quad \delta^*(P_{c})=1.
\end{split}
\end{equation*}

\item {\sl Octahedron}, $4\,|\,2\,3$, $f=(6,12,8)$ (cf.~\cite{m}).
\begin{equation*}
\begin{split} 
&  P_{o}      = \left\{x \in \R^3 :  |x_1|+|x_2|+|x_3| \leq 1 \right\},\\
&  \Lambda^*(P_{o})  =  2\left( (\frac{1}{3}, \frac{1}{2}, \frac{1}{6})^\intercal,
                (-\frac{1}{6}, -\frac{1}{3}, \frac{1}{2})^\intercal,
                (-\frac{1}{2}, \frac{1}{6}, -\frac{1}{3})^\intercal 
                 \right)\Z^3, \\
&  \delta^*(P_{o})  = \frac{18}{19}\approx 0.947368421.
\end{split}
\end{equation*}

\item {\sl Dodecahedron}, $3\,|\,2\,5$, $f=(20,30,12)$.
 \begin{equation*}
\begin{split} 
  & P_{d}    = \left\{x \in \R^3 : |\tau x_1| + |x_{2}| \leq 1,\, 
                                   |\tau x_2| + |x_{3}| \leq 1,\, 
                                   |\tau x_3| + |x_{1}| \leq 1 \right\},  \\
& \Lambda^*((1+\tau)P_{d})  =  2D_3,\quad \delta^*(P_{d})   = \frac{2+\tau}{4}\approx 0.904508497.
\end{split}
\end{equation*}

\item {\sl Icosahedron}, $5\,|\,2\,3$, $f=(12,30,20)$.
\begin{equation*}
\begin{split}
& P_{i}  = \big\{x\in \R^3 : |x_1| + |x_2| + |x_3| \leq 1, \,
                             |\tau x_1| + |(1/\tau)x_3| \leq 1,\, \\
&\quad\quad\quad\quad\quad\quad\quad   |\tau x_2| + |(1/\tau)x_{1}| \leq 1,\, 
                             |\tau x_3| + |(1/\tau)x_{2}| \leq 1  \big\},\\
& \Lambda^*((1+\tau)P_{i})  =  2(w^1({\bar x}),w^2({\bar x}),
w^3({\bar x}))\Z^3 \text{ where }
\end{split}
\end{equation*}
\vspace{-0.2cm}
\begin{equation}
\label{latt_icosi}
\begin{split}
& w^1({\bar x})=\begin{pmatrix}
( - { \frac {33}{8}}  - { \frac {39}{8}
} \,\sqrt{5})\,{\bar x}^{2} + ({ \frac {39}{4}}  + 
{ \frac {33}{4}} \,\sqrt{5})\,{\bar x} - { 
\frac {11}{4}}  - { \frac {3}{2}} \,\sqrt{5} \\
( - { \frac {1}{4}}  - { \frac {1}{4}} 
\,\sqrt{5})\,{\bar x} + 1 + { \frac {1}{2}} \,\sqrt{5} \\
({ \frac {33}{8}}  + { \frac {39}{8}} 
\,\sqrt{5})\,{\bar x}^{2} + ( - { \frac {19}{2}}  - 8\,
\sqrt{5})\,{\bar x} + { \frac {13}{4}}  + { 
\frac {3}{2}} \,\sqrt{5}
\end{pmatrix}, 
\\[0.5ex]
& w^2({\bar x})=\begin{pmatrix}
( - { \frac {33}{40}} \,\sqrt{5} - { 
\frac {39}{8}} )\,{\bar x}^{2} + ({ \frac {41}{20}} \,
\sqrt{5} + { \frac {35}{4}} )\,{\bar x} - { 
\frac {5}{2}}  - { \frac {23}{20}} \,\sqrt{5} \\
({ \frac {5}{4}}  + { \frac {1}{4}} \,
\sqrt{5})\,{\bar x} - 1 - { \frac {1}{2}} \,\sqrt{5} \\
( - { \frac {33}{40}} \,\sqrt{5} - { 
\frac {39}{8}} )\,{\bar x}^{2} + ({ \frac {9}{5}} \,\sqrt{5
} + { \frac {15}{2}} )\,{\bar x} - { \frac {3
}{20}} \,\sqrt{5}
\end{pmatrix}, 
\\[0.5ex]
& w^3({\bar x})=
\begin{pmatrix}
({ \frac {1}{2}} \,\sqrt{5} + { \frac {
3}{2}} )\,{\bar x} - 2 - \sqrt{5}, &
{\bar x}, &
0
\end{pmatrix}^\intercal, 
\end{split}
\end{equation}
\vspace{-0.2cm}
\begin{equation*}
\begin{split}
& \text{and } {\bar x} \text{ is the unique root with } {\bar x}\in(1,2) \text{ of } \\
& 1086\,x^3+(-1063-111\sqrt{5})\,x^2+(15\sqrt{5}+43)\,x+102+44\sqrt{5}, \\
& \delta^*(P_{i})=\frac{5(1+\tau)}{|\det(w^1({\bar x}),w^2({\bar x}),w^3({\bar x}))|}\approx 0.836357445.
\end{split}
\end{equation*}

\item {\sl Cubeoctahedron}, $2\,|\,3\,4$, $f=(12,24,14)$ (cf.~\cite{h}).
\begin{equation*}
\begin{split}
&P_{co}=\left\{x\in \R^3 : x \in P_{c} \cap 2\cdot P_{o}\right\}, \\
&\Lambda^*(P_{co})=\Lambda^*(P_{t}),\quad \delta^*(P_{co})=\frac{45}{49}\approx 0.918367346.
\end{split}
\end{equation*}

\item {\sl Icosidodecahedron}, $2\,|\,3\,5$, $f=(30,60,32)$.
\begin{equation*}
\begin{split}
&P_{id}=\left\{x\in \R^3 : x \in P_{i} \cap  P_{d}\right\}, \\
&\Lambda^*((1+\tau)P_{id})=2D_3,\quad \delta^*(P_{id})=\frac{14+17\tau}{48}\approx 0.864720371.
\end{split}
\end{equation*}

\item  {\sl Rhombic Cubeoctahedron}, $3\,4\,|\,2$, $f=(24,48,26)$.
\begin{equation*}
\begin{split}
&P_{rco}=\big\{x\in \R^3 : |x_1|+|x_2|\leq 2,\, |x_1|+|x_3|\leq 2,\, 
                              |x_2|+|x_3|\leq 2, \\
&\quad\quad\quad\quad\quad\quad\quad\quad\text{ and } x\in \sqrt{2}\cdot P_{c}\cap 
                                        (4-\sqrt{2})P_{o}\big\},\\
&\Lambda^*(P_{rco})=2D_3,\quad \delta^*(P_{rco})=\frac{16\sqrt{2}-20}{3}\approx 0.875805666.
\end{split}
\end{equation*}

\item {\sl Rhombic Icosidodecahedron}, $3\,5\,|\,2$, $f=(60,120,62)$.
\begin{equation*}
P_{rid}=\Big\{ x\in \R^3 : x\in (3\tau+2)\cdot P_{rtc} \cap (4\tau+1)\cdot P_{i} \cap (3(1+\tau))\cdot P_{d} \Big\}, 
\end{equation*}
\vspace{-1.0truecm}
\begin{equation*}
\begin{split}
&\Lambda^*(P_{rid})=2\Big( \left(\frac{\tau-1}{4\tau+2}, \frac{7}{2}, \frac{9\tau+4}{4\tau+2}\right)^\intercal, \left(\frac{9\tau+4}{4\tau+2}, \frac{\tau-1}{4\tau+2}, \frac{7}{2} \right)^\intercal, \\
&\hspace{6.2truecm}\left(\frac{7}{2}, \frac{9\tau+4}{4\tau+2}, \frac{\tau-1}{4\tau+2} \right)^\intercal\Big)\Z^3,\\ 
&\delta^*(P_{rid})=\frac{8\tau+46}{36\tau+15}\approx 0.804708487.
\end{split}
\end{equation*}

\item {\sl Truncated Cube}, $2\,3\,|\,4$, $f=(24,36,14)$.
\begin{equation*}
\begin{split}
&P_{trc}=\left\{ x\in \R^3: x \in P_{c}\cap (1+\sqrt{2})P_{o}\right\}, \\
&\Lambda^*(P_{trc})=2\left( (1, -\alpha, 0)^\intercal,
                (0,1,-\alpha)^\intercal,
                (-\alpha,0,1)^\intercal
                 \right)\Z^3, \quad  
 \alpha=(2-\sqrt{2})/3, \\
&\delta^*(P_{trc})=\frac{9}{5+3\sqrt{2}}\approx 0.973747688.
\end{split}
\end{equation*}

\item {\sl Truncated Octahedron}, $2\,4\,|\,2$, $f=(24,36,14)$.
\begin{equation*}
\begin{split}
&P_{tro}=\left\{ x \in \R^3 : x \in P_{c} \cap \frac{3}{2} P_{o} \right\}, \\
&\Lambda^*(P_{tro})=2\left( (1,0,0)^\intercal, (1, 1, 0)^\intercal, (\frac{1}{2}, \frac{1}{2},-\frac{1}{2})^\intercal\right)\Z^3,\quad \delta^*(P_{tro})=1.
\end{split}
\end{equation*}

\item {\sl Truncated Dodecahedron}, $2\,3\,|\,5$, $f=(60,90,32)$.
\begin{equation*}
\begin{split}
&P_{trd}=\{ x\in \R^3 : x \in (1+\tau)\cdot P_d \cap \left((7+12\tau)/(3+4\tau)\right)\cdot P_i\}, \\
&\Lambda^*(P_{trd})=2D_3, \quad 
\delta^*(P_{trd})=\frac{1}{4}\cdot \frac{5\tau+16}{6\tau-3}\approx 0.897787626.
\end{split}
\end{equation*}

\item {\sl Truncated Icosahedron}, $2\,5\,|\,3$, $f=(60,90,32)$.
\begin{equation*}
\begin{split}
&P_{tri}=\{ x\in \R^3 : x \in (1+\tau)\cdot P_i \cap (4/3+\tau)\cdot P_d\}, \\
&\Lambda^*(P_{tri})=\Lambda^*((1+\tau)P_{i}), \\
&\delta^*(P_{tri})=\frac{43\sqrt{5}+125}{108\cdot |\det(w^1({\bar x}),w^2({\bar x}),w^3({\bar x}))|} \approx 0.7849877759.\quad (\text{cf.~\eqref{latt_icosi}})
\end{split}
\end{equation*}

\item {\sl Truncated Cubeoctahedron}, $2\,3\,4\,|$, $f=(48,72,26)$.
\begin{equation*}
\begin{split}
&P_{trco}=\big\{x\in \R^3 : |x_1|+|x_2| \leq 2+3\sqrt{2},\, |x_2|+|x_3| \leq 2+3\sqrt{2}, \\ 
&\quad\quad\quad |x_2|+|x_3| \leq 2+3\sqrt{2} \,\text{ and } 
 x\in (2\sqrt{2}+1)P_c\cap (3\sqrt{2}+3)P_0\big\},  \\
& \Lambda^*(P_{trco})=2\Big( (2\sqrt{2}+1, -2\sqrt{2}-\frac{1}{2}+\alpha, 2\sqrt{2}+\frac{1}{2}-\alpha)^\intercal, \\
&\quad\quad\quad (\frac{1}{4}\sqrt{2}-\frac{3}{4}+\frac{1}{2}\alpha, -\frac{3}{4}\sqrt{2}+\frac{1}{4}+\frac{1}{2}\alpha, 2\sqrt{2}+1)^\intercal, \\
&\quad\quad\quad (\frac{7}{4}+\frac{7}{4}\sqrt{2}-\frac{1}{2}\alpha, \frac{1}{2}+\alpha, 
\frac{5}{4}\sqrt{2}+\frac{3}{4}-\frac{1}{2}\alpha)^\intercal\Big)\Z^3, \\
& \text{ where } \alpha = \frac{1}{6}\sqrt{33}(\sqrt{2}+1), \\
&\delta^*(P_{trco})=\frac{99}{992}*\sqrt{66}-\frac{231}{1984}*\sqrt{33}+
\frac{2835}{992}*\sqrt{2}-\frac{6615}{1984} 
\approx 0.849373252.
\end{split}
\end{equation*}

\item {\sl Truncated Icosidodecahedron}, $2\,3\,5\,|$, $f=(120,180,62)$.
\begin{equation*}
\begin{split}
&P_{trid}=\Big\{ x\in \R^3 : x\in (5\tau+4)\cdot P_{rtc}\cap (6\tau+3)\cdot P_{i} \cap (5(1+\tau))\cdot P_{d} \Big\} \\
&\Lambda^*(\frac{1}{5} P_{trid})=2D_3, \quad
\delta^*(P_{trid})=\frac{2}{5}\tau+\frac{9}{50}\approx 0.827213595.
\end{split}
\end{equation*}

\item {\sl Truncated Tetrahedron}, $2\,3\,|\,3$, $f=(12,18,8)$.
\begin{equation*}
\begin{split}
&P_{trt}=\{ x\in \R^3 : x\in 5\cdot P_t \cap -3\cdot P_t \},\\
&\Lambda^*(P_{trt})=2\left( (\frac{2}{3},2,\frac{4}{3})^\intercal, 
        (2,-\frac{4}{3}-\frac{2}{3})^\intercal,
(-\frac{4}{3},\frac{2}{3},-2)^\intercal\right)\Z^3, \\
&\delta^*(P_{trt})=\frac{207}{304}\approx 0.680921053.
\end{split}
\end{equation*}

\item {\sl Snub Cube}, $|\,2\,3\,4$, $f=(24,60,38)$.
\par Let $P_{sc}$ be the snub cube such that the 6 quadrangle facets lie in the hyperplanes 
$\{x\in \R^3 : x_i=\pm 1\}$, $1\leq i \leq 3$, and let $y^*$ be the unique real solution of $y^3+y^2+y=1$.  
\begin{equation*}
\begin{split}
\Lambda^*(P_{sc})&=2\cdot\left( (1,0,0)^\intercal, (0,0,1)^\intercal, (\frac{1}{2},\frac{1}{y^*}-1,-\frac{1}{2})^\intercal\right),\\
\delta^*(P_{sc})&=\frac{1}{2}+\frac{1}{6}y^*+\frac{2}{3}(y^*)^2\approx 0.78769996.
\end{split}
\end{equation*}

\item {\sl Snub Dodecahedron}, $|\,2\,3\,5$, $f=(60,150,92)$.
\par Let $P_{sd}$ be the snub dodecahedron such that the 12 pentagonal facets lie in the supporting hyperplanes of the facets of the dodecahedron $(1+\tau)P_d$. 
\begin{equation*}
\begin{split}
 \Lambda^*(P_{sd})=2D_3,\quad \delta^*(P_{sd})={\rm vol}(P_{sd})/16\approx 0.788640117.
\end{split}
\end{equation*}
\end{itemize}
\vspace{0.5cm}
\begin{theorem} 
\label{final_theorem}
The above list contains the densities of a densest lattice packing of all regular and Archimedean polytopes.
\end{theorem}
\begin{proof} Algorithm \ref{algorithm_final}.
\end{proof}

\noindent{\bf Remarks.} 
\begin{enumerate}
\item[a)]
In order to get exact values for the densities as given in the above list we first use a numerical implementation of Algorithm \ref{algorithm_final}, by which we determine an  optimal selection $F_{l_1},\dots,F_{l_k}$ of facets corresponding to a critical lattice. Then for this special choice of factes we carry out the steps R1., R2.~and R3.~with the symbolic computer algebra system {\tt Maple V Release 5}.   
\item[b)]
 The last three polytopes $P_{trt}, P_{sc}, P_{sd}$ are not centrally symmetric and thus one has to calculate the difference bodies first. The difference body of $P_{sc}$ is a polytope with 74 facets, whereas $\frac{1}{2}(P_{sd}-P_{sd})$ has already 182 facets.
 We also have computed the densest packing lattice of the dual polytopes of the Archimedean polytopes, the so called Catalan polytopes. Thereby we had to determine packing lattices of polytopes with more than 380 facets. The CPU-time for the determination of the densest packing lattices of all regular and Archimedean polytopes is about 5.5 hours on a PC with a 266\,Mhz-Pentium II processor.  
\end{enumerate}

\ifnum\pictures=1
\vspace{0.0cm}
\noindent 
On the following pages we present pictures of the regular and Archimedean polytopes. The pictures were generated with the program {\tt povray} (see http://www.povray.org). The dark points are the lattice points of a critical lattice of the difference body lying on the boundary of this body.

\vspace{0.8cm}
\begin{minipage}[h]{.46\linewidth}
\centering\epsfig{file=\wherepoly/tetrahedron.ps, width=\linewidth}
\centering{Tetrahedron}
\hphantom{(Difference body of Tetrahedron)} 
\end{minipage}\hfill
\begin{minipage}[h]{.46\linewidth}
\centering\epsfig{file=\wherepoly/cubeoctahedron_pack.ps, width=\linewidth}
\centering{Difference body of Tetrahedron (Cubeoctahedron)} 
\end{minipage}

\vspace{0.5cm}
\begin{minipage}[h]{.46\linewidth}
\centering\epsfig{file=\wherepoly/cube_pack.ps, width=\linewidth}
\centering{Cube} 
\end{minipage}\hfill
\begin{minipage}[h]{.46\linewidth}
\centering\epsfig{file=\wherepoly/crosspolytope_pack.ps, width=\linewidth}
\centering{Octahedron} 
\end{minipage}

\vspace{0.2cm}
\begin{minipage}[h]{.46\linewidth}
\centering\epsfig{file=\wherepoly/dodecahedron_pack.ps, width=\linewidth}
\centering{Dodecahedron} 
\end{minipage}\hfill
\begin{minipage}[h]{.46\linewidth}
\centering\epsfig{file=\wherepoly/icosahedron_pack.ps, width=\linewidth}
\centering{Icosahedron} 
\end{minipage}

\vspace{0.5cm}
\begin{minipage}[h]{.46\linewidth}
\centering\epsfig{file=\wherepoly/cubeoctahedron_pack.ps, width=\linewidth}
\centering{Cubeoctahedron} 
\end{minipage}
\hfill
\begin{minipage}[h]{.46\linewidth}
\centering\epsfig{file=\wherepoly/icosidodecahedron_pack.ps, width=\linewidth}
\centering{Icosidodecahedron} 
\end{minipage}

\vspace{0.3cm}
\begin{minipage}[h]{.46\linewidth}
\centering\epsfig{file=\wherepoly/rhombicubeoctahedron_pack.ps, width=\linewidth}
\centering{Rhombic Cubeoctahedron} 
\end{minipage}
\hfill
\begin{minipage}[h]{.46\linewidth}
\centering\epsfig{file=\wherepoly/rhombiicosidodecahedron_pack.ps, width=\linewidth}
\centering{Rhombi Icosidodecahedron} 
\end{minipage}

\vspace{0.5cm}
\begin{minipage}[h]{.46\linewidth}
\centering\epsfig{file=\wherepoly/truncated_cube_pack.ps, width=\linewidth}
\centering{Truncated Cube} 
\end{minipage}
\hfill
\begin{minipage}[h]{.46\linewidth}
\centering\epsfig{file=\wherepoly/truncated_crosspolytope_pack.ps, width=\linewidth}
\centering{Truncated Octahedron} 
\end{minipage}

\vspace{0.2cm}
\begin{minipage}[h]{.46\linewidth}
\centering\epsfig{file=\wherepoly/truncated_dodecahedron_pack.ps, width=\linewidth}
\centering{Truncated Dodecahedron} 
\end{minipage}
\hfill
\begin{minipage}[h]{.46\linewidth}
\centering\epsfig{file=\wherepoly/truncated_icosahedron_pack.ps, width=\linewidth}
\centering{Truncated Icosahedron} 
\end{minipage}

\vspace{0.3cm}
\begin{minipage}[h]{.46\linewidth}
\centering\epsfig{file=\wherepoly/truncated_cubeoctahedron_pack.ps, width=\linewidth}
\centering{Truncated Cubeoctahedron} 
\end{minipage}
\hfill
\begin{minipage}[h]{.46\linewidth}
\centering\epsfig{file=\wherepoly/truncated_icosidodecahedron_pack.ps, width=\linewidth}
\centering{Truncated Icosidodecahedron} 
\end{minipage}

\vspace{0.3cm}
\begin{minipage}[h]{.46\linewidth}
\centering\epsfig{file=\wherepoly/truncated_tetrahedron.ps, width=\linewidth}
\centering{Truncated Tetrahedron}
\hphantom{hugohugohugo} 
\end{minipage}
\hfill
\vspace{0.5cm}
\begin{minipage}[h]{.46\linewidth}
\centering\epsfig{file=\wherepoly/truncated_tetsym_pack.ps, width=\linewidth}
\centering{Difference body of Truncated Tetrahedron} 
\end{minipage}

\begin{minipage}[h]{.46\linewidth}
\centering\epsfig{file=\wherepoly/snub_cube.ps, width=\linewidth}
\centering{Snub Cube} 
\end{minipage}
\hfill
\begin{minipage}[h]{.46\linewidth}
\centering\epsfig{file=\wherepoly/snub_cube_diff.ps, width=\linewidth}
\centering{Difference body of Snub Cube} 
\end{minipage}

\vspace{0.5cm}
\begin{minipage}[h]{.46\linewidth}
\centering\epsfig{file=\wherepoly/snub_dodecahedron.ps, width=\linewidth}
\centering{Snub Dodecahedron} 
\hphantom{hugohugohugo}
\end{minipage}
\hfill
\begin{minipage}[h]{.46\linewidth}
\centering\epsfig{file=\wherepoly/snub_dodecahedron_diff.ps, width=\linewidth}
\centering{Difference body of Snub Dodecahedron} 
\end{minipage}

\fi

\vspace{0.4cm}
\noindent
{\it Acknowledgements.} We would like to thank Utz-Uwe Haus, Achill Sch\"ur\-mann and  Chuanming Zong for helpful comments.

\vspace{0.4cm}
\noindent
{\sc \small 
Ulrich Betke, University of Siegen, Department of Mathematics,   \linebreak D-57068 Siegen, Germany.} 

\noindent {\it E-mail:} {\tt betke@mathematik.uni-siegen.de}

\vspace{0.2cm}
\noindent
{\sc \small Martin Henk, University of Magdeburg, Department of Mathematics/IMO,  D-39106 Magdeburg, Germany.} 

\noindent{\it E-mail:} {\tt henk@imo.math.uni-magdeburg.de} 



\end{document}